\documentclass{autart}
\usepackage{comment}
\usepackage{amsmath,amssymb,amsfonts}
\usepackage{algorithmic}
\usepackage{graphicx}
\usepackage{textcomp}
\usepackage{mathrsfs}
\usepackage{bm}
\usepackage{xcolor}
\usepackage{ntheorem}
\newtheorem*{Proof}{Proof}       
\newtheorem{Proposition}{\textbf{Proposition}} 
\newtheorem{Notation}{\textbf{Notation}} 

\newtheorem{Lemma}{\textbf{Lemma}}
\newtheorem{theorem}{\textbf{Theorem}} 
\newtheorem{Remark}{\textbf{Remark}} 
\newtheorem{Assumption}{\textbf{Assumption}} 
\allowdisplaybreaks[4]
\usepackage{epsfig}
\usepackage{subfig}

\usepackage[title]{appendix}
\usepackage{float} 
\graphicspath{{pictures/}}
\begin{document}
\setlength\abovedisplayskip{6.4pt}
        \setlength\belowdisplayskip{6.4pt}
        \setlength\abovedisplayshortskip{6.4pt}
        \setlength\belowdisplayshortskip{6.4pt}
        \allowdisplaybreaks
        \setlength{\parindent}{1em}
        \setlength{\parskip}{-0em}  
        \addtolength{\oddsidemargin}{6.4pt} 
\begin{frontmatter}

\title{Robustness of Reaction-Diffusion PDEs Predictor-Feedback to Stochastic Delay Perturbations \thanksref{t2}} 

\thanks[t2]{Dandan Guan and Jie Qi are supported by National Natural Science Foundation of China (62173084, 61773112) and the Project of Science and Technology Commission of Shanghai Municipality, China (23ZR1401800, 22JC1401403), in part by the Fundamental Research Funds for the Central Universities and Graduate Student Innovation Fund of Donghua University CUSF-DH-D-2023043. Mamadou Diagne is supported by National Science
Foundation CAREER Award CMMI-2302030. This work was performed during Dandan Guan's visit to the University of California  San Diego.}

\thanks[footnoteinfo]{Corresponding author M. Diagne Tel. +1(518)8925350.}

\author[Paestum]{Dandan Guan}\ead{gdd\_dhu@163.com},   
\author[Paestum]{Jie Qi}\ead{ jieqi@dhu.edu.cn}\ and  
\author[SD]{Mamadou Diagne \thanksref{footnoteinfo}}\ead{mdiagne@ucsd.edu}

\address[Paestum]{College of Information Science and
	Technology, Donghua University, Shanghai, China, 201620} 
 \address[SD]{Department of Mechanical and Aerospace Engineering, University of California San Diego, La Jolla, CA, USA, 92093}

\begin{abstract}
This paper studies the robustness of a PDE backstepping delay-compensated boundary controller for a reaction-diffusion partial differential equation (PDE)  with respect to a nominal delay subject to stochastic  error disturbance.    The stabilization problem under consideration involves random perturbations  modeled by a finite-state Markov process that further  obstruct the actuation path at the controlled boundary of the infinite-dimension plant. This scenario is useful to describe several actuation failure modes in process control. Employing the recently introduced infinite-dimensional representation of the state of an actuator subject to stochastic input delay for ODEs (Ordinary Differential Equations), we convert the stochastic input delay into  $r+1$  unidirectional advection PDEs, where $r$ corresponds to the number of jump states.  Our stability analysis assumes full-state measurement of the spatially distributed plant's state and relies on a hyperbolic-parabolic PDE cascade representation of the plant plus actuator dynamics. Integrating the plant and the nominal stabilizing boundary control action, all while considering probabilistic delay disturbances, we establish the proof of mean-square exponential stability 
 as well as the well-posedness of the closed-loop system when  random phenomena weaken the nominal  actuator compensating effect.   Our proof is based on the Lyapunov method, the theory of infinitesimal operator for stability, and $C_0$-semigroup  theory for well-posedness. Our stability result refers to the $L^2$-norm of the plant state and the $H^2$-norm of the actuator state. Furthermore, our study presents a qualitative analysis of the maximum deviation  of the stochastic disturbance relative to the nominal delay, which is expressed as a function of the severity of the plant instability--the magnitude of the \emph{reactivity coefficient}-- and the known \emph{upper bound of the transition rate} of the underlying stochastic perturbations.  Extensive simulation results illustrate the viability of the proposed robustness study.
\end{abstract}

\end{frontmatter}
\section{Introduction}
\subsection{Stabilization of parabolic PDE systems with delays}
The omnipresence of time-delay in engineering systems modeled by partial differential equations (PDEs) or delay differential equations (DDEs) poses both theoretical and practical challenges that are yet to be overcome by control engineers.  PDE backstepping is a powerful tool that has been widely used to design predictor-feedback controllers for linear and nonlinear finite-dimensional systems with various types of input delays. It has proven to be effective in cases where the delays are known (constant, time-dependent, delayed input-dependent or state-dependent) \cite{krstic2009input,bekiaris2013,diagne2017compensation,diagne2017time,diagne2017control} and cases where the delays are unknown \cite{bresch2009}. An initial design of backstepping to infinite-dimensional systems with delays was the design of a predictor-feedback controller that compensated for a boundary input delay of a scalar 1D reaction-diffusion system \cite{Krstic2009control}. Indeed, \cite{Krstic2009control} achieves exponential stabilization of both the plant and actuator states by introducing a hyperbolic-parabolic PDE-PDE cascade representation.
 Subsequently, the result was further expanded to diffusive 3D continuum formation control under cylindrical topology \cite{qi2019control}. The approach proposed in \cite{Krstic2009control} has resulted in meaningful progress towards the development of stabilizing controllers for delay-dependent hyperbolic and parabolic partial differential equations (PDEs), crowned with many interesting contributions \cite{deutscher2019fredholm,qi2020compensation,Wang2021adaptive}. As a result of \cite{qi2020compensation}, a feedback compensator for a reaction-diffusion PDE with spatially varying delays affecting its distributed input was developed, while \cite{Wang2021adaptive,Wang2022,wang2023observer} designed delay-adaptive boundary and in-domain stabilizing control laws  despite the presence of an arbitrarily large, constant and unknown delay. An alternative approach relying on LMI conditions  for achieving exponential stability of reaction-diffusion PDEs by combining a finite-dimensional observer with a predictor was introduced in \cite{katz2021sub,katz2020constructive}.  Initially  presented in \cite{katz2020constructive}, the method was generalized in \cite{lhachemi2022predictor} to allow for the stabilization of systems subject to Neumann boundary control and measurement.

\subsection{Control of ODE systems with probabilistic  delays}

Several contributions have confirmed the occurrence of probabilistic delays in power systems \cite{lu2015mathematical,liu2018stability} and traffic flow control \cite{orosz2006traffic,sykora2020moment}. These studies are often  motivated by the probabilistic nature of delays caused by  random packets arrival in networked control systems  \cite{krtolica1994stability,Liang2010,hu2017robust,kaviarasan2016robust} including biological gene regulatory networks \cite{gomez2016stability,lou2010exponential}. Results pertaining to the design of robust control laws for linear systems affected by a random delay are available in \cite{krtolica1994stability,chen2003guaranteed,yue2009stabilization}. In these works, the delay is defined as a known probability distribution rather than being considered as a constant  mismatch \cite{mondie2003finite,mondie2001delay,krstic2008lyapunov} or as any other class of uncertainties \cite{kim1996robust,choi1996observer,bekiaris2013robustness,karafyllis2013delay}.
  Employing 
  a  structural perturbation approach, \cite{krtolica1994stability}  deduced the exponential stability of linear discrete-time systems with Markov chain delays.        In \cite{gomez2016stability},  stability conditions of a nonlinear plant undergoing stochastic delays are derived from the mean and the second moment dynamics  (see \cite{sadeghpour2019stability} and references therein).  Inspired by the probability delay averaging method introduced in \cite{kolmanovsky2001mean} and \cite{kolmanovsky2001optimal}, the authors of \cite{kong2022prediction1} and \cite{kong2022prediction2} recently proposed robust compensators for both linear and nonlinear systems with a nominal input delay  subject to a random disturbance.  Assuming that the stochastic perturbation values remain relatively small and closely aligned with the prediction horizon, \cite{kong2022prediction1,kong2022prediction2} establish that employing PDE backstepping predictor-feedback designed for constant delay  \cite{krstic2008lyapunov}, is adequate for ensuring mean-square exponential stability for both the plant's and actuator state. The stochastic delay representation in \cite{kong2022prediction1,kong2022prediction2}  is relevant when the actuator state follows a Markov process with a finite state space. The method allows to capture the effect  of the continuous-time Markov jump process  by breaking down the system under consideration into a nominal delay used for constructing the nominal compensator and an additional stochastic delay error disturbance generated by the randomly varying perturbation.
\subsection{Contributions}
 Despite all the progress made on the design of predictor-feedback control laws for finite- and infinite-dimensional systems, the design of controllers for PDE systems with stochastic delays or the study of the robustness of nominal delay compensators to stochastic perturbations has remained an open challenge. Borrowing the actuator's state representation introduced in \cite{kong2022prediction1} for ODEs with stochastic input delays, we study the robustness of the controller \cite{Krstic2009control} when a stochastic error disturbance described by a Markov jump process affects the boundary input delay. The closed-loop system subject to relatively small stochastic perturbations achieves mean-square exponential stability with respect to the $L^2$-norm of the plant's state and $H^2$-norm of the actuator state. Our proof relies on the construction of a target system whose mean-square exponential stability is established with the help of  Lyapunov arguments and infinitesimal operator.  Our proof of exponential stabilization is valid when the expectation of the time-varying delay mismatch driven by stochastic transitions complies with a constraint linked to the plant's instability severity. In other words, it depends on the magnitude of the coefficient of reactivity and the upper limit for the transition rate.
 The proof of the well-posedness of the closed-loop system, which holds regardless of the infrequency of jumps in the Markov process or the assigned maximum transition rate, is provided using  $C_0$-semigroup theory.  
The presented stability and well-posedness analysis are not directly comprehensible by straightforwardly borrowing ODEs'  state feedback tools \cite{kong2022prediction1,kong2022prediction2}. Our stability analysis differs from \cite{Krstic2009control}  and the proof of the well-posedness of the closed-loop system does not result from a trivial analysis in that random factors are considered.

This paper is organized as follows. Section \ref{2} presents the design of the boundary controller for reaction-diffusion PDEs with stochastic delay. Section \ref{6} studies the well-posedness of the system by semigroup theory. Section \ref{3} states the stability analysis and the numerical simulations to support the theoretical statements are provided in Section \ref{4}. The paper ends with concluding remarks and future perspectives in Section \ref{5}.
\begin{Notation}\rm
	In this paper, we define the norms for $f(x)\in L^2(0,1)$ and  $\bm{v}(x)=(v_1(x),v_2(x),\cdots,v_r(x))^T$, where $v_i(x)\in L^2(0,1)$, $i=1,\cdots,r$, $r\in\mathbb{N}^+$,
	\begin{align}
	\Vert f\Vert_{L^2}^2&=\int_{0}^{1}f^2(x)dx\\
	\Vert f\Vert_{H^1}^2&=\Vert f\Vert_{L^2}^2+\Vert f_x\Vert_{L^2}^2\\
	\Vert \bm{v}\Vert_{L^2}^2&=\int_{0}^{1}\bm{v}^T(x)\bm{v}(x)dx=\sum_{i=1}^{r}\int_{0}^{1}v^2_i(x)dx\\
	\Vert \bm{v}\Vert_{H^1}^2&=\Vert \bm{v}\Vert_{L^2}^2+\Vert \bm{v}_x\Vert_{L^2}^2.
	\end{align} 
	
\end{Notation}

\section{Backstepping Boundary Control Design}\label{2}
\subsection{Problem statement }
We consider  the following potentially unstable reaction-diffusion  PDE with a selected nominal delay subject to  stochastic boundary input delay
\begin{align} 
\label{u-main}
u_{t}(x,t)&=u_{xx}(x,t)+\lambda u(x,t)\\
\label{u-bnd0}
u(0,t)&=0\\
\label{u-control}
u(1,t)&=U(t-\tau(t))\\
u(x,0)&=u_0(x)\label{initial-u}
\end{align}
where $x\in [0,1]$, $t>0$, $\lambda>0$.  Here, $U(t)$ is the boundary control input subject to a  stochastic delay $\tau(t),$ which 
  is described  by a Markov process  with a transition probability $P_{ij}(s,t)$, $t>s$, that is to say, the probability to switch from $\tau_i$ at time $s$ to $\tau_j$ at time $t$ satisfies the Kolmogorov forward differential equation (see. pages 19--23 in Reference \cite{yin2012continuous} and Reference \cite{zeifman2019study}). That is
\begin{align}\label{p-pi}
\frac{\partial P_{ij}(s,t)}{\partial t}&=\sum_{k=1}^{r}q_{kj}(t)P_{ik}(s,t),~i,j=1,2,\cdots,r\\
P_{ii}(s,s)&=1,~P_{ij}(s,s)=0,~j\ne i
\end{align}
where the transition rate $q_{ij}(t)$  in \eqref{p-pi} satisfies
\begin{align}
    -q_{ii}(t)=&\lim_{\Delta t\rightarrow 0}\frac{1-P_{ii}(t,t+\Delta t)}{\Delta t}\\
    q_{ij}(t)=&\lim_{\Delta t\rightarrow 0}\frac{P_{ij}(t,t+\Delta t)}{\Delta t},~j\ne i
    \\
q_{ii}(t)=&-\sum_{j=1,j\ne i}^{r}q_{ij}(t).\label{10}
\end{align} 
In addition, the Markov process $\tau(t)$ has the right continuous realizations with $\tau(t)\in\{\tau_i|0<\underline{\tau}\le \tau_i<\tau_j\le\overline{\tau}, i <j; i,j\in\{1,\cdots,r\}\}$ having a finite-dimensional state. \\
\begin{Assumption}[\cite{yin2012continuous}]\label{assumption1}\rm
		  The transition rate $ q_{ij}(t)$ satisfies $\vert q_{ij}(t)\vert\le q^\star$, $i,j\in\{1,2,\cdots,r\}$, where  $q^\star$ is a known positive constant.
\end{Assumption}
Assumption \ref{assumption1} is realistic  as a known upper bound $q^\star$ of the transition rate is often available \cite{yin2012continuous}.\\
\begin{Assumption}\label{assumption2}\rm
		  Considering  a known nominal delay $\tau_0\in[\underline{\tau},\overline{\tau}]$,  it is assumed that the stochastic delay satisfies 
    \begin{align}
        \mathbb{E}_{[0,\tau(0)]}\left(|\tau_0-\tau(t)|^2\right)\le\varepsilon
    \end{align}
    where $\varepsilon$ is a sufficiently small value, $\mathbb{E}_{[0,\tau(0)]}\left(|\tau_0-\tau(t)|^2\right)$ is the condition expectation of $|\tau_0-\tau(t)|^2$ at instant $t$ under the condition $\tau(t)=\tau(0)$ at $t=0$.
\end{Assumption}
 The precise characterization of $\varepsilon$  is given in the stability analysis.
 
To design a backstepping controller,  let us first define the following vector functions \cite{kong2022prediction1,kong2022prediction2} 
\begin{align}
\bm{v}(x,t)=&(v_1(x,t), v_2(x,t),\cdots, v_r(x,t))^T\label{act}\\
\hat{v}(x,t)=&U(t+\tau_0(x-1))\label{constant-act}\\
\bm{\tilde{v}}(x,t)=&\bm{v}(x,t)-\bm{1}\hat{v}(x,t)\label{error-control}
\end{align}
for all $ t \in [t-\tau_0,t]$, where  $\tau_0\in[\underline{\tau},\overline{\tau}]$ is a constant delay, $v_i(x,t)=U(t+\tau_i(x-1))$, $i=1,2,\cdots,r$, and $\bm{1}$ is the r-by-1 all-ones vector.

\begin{figure*}[htbp]
	\centering
	\centering{\epsfig{figure=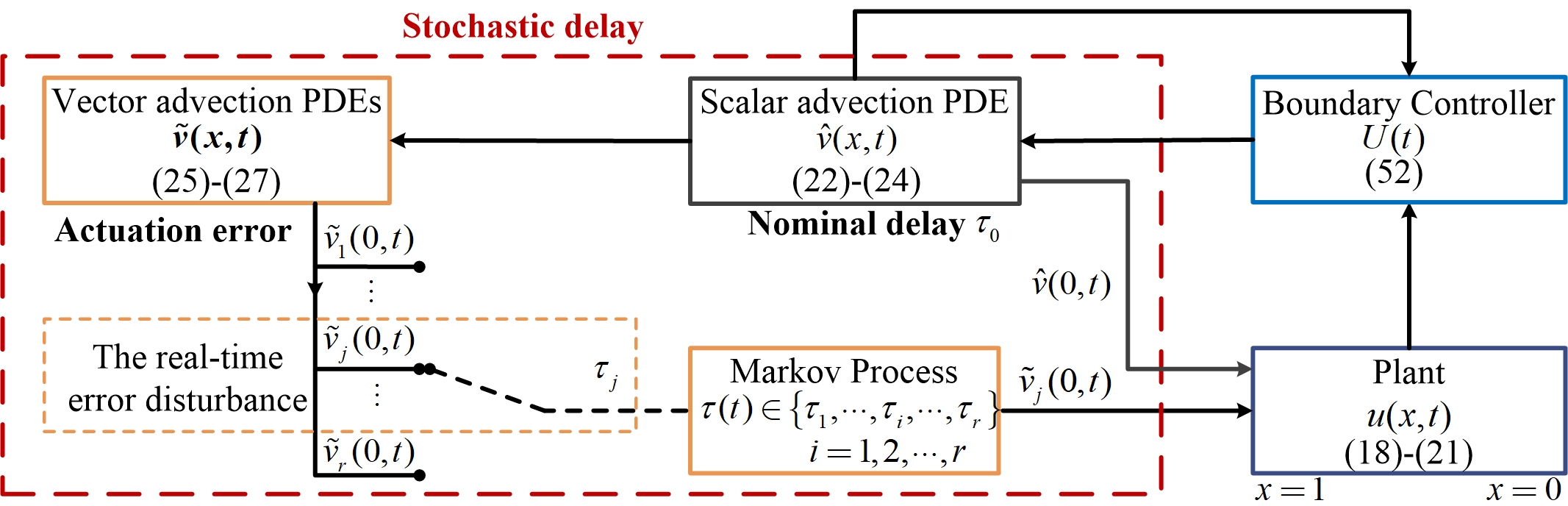,width=0.85\textwidth}}
	\caption{The structure of the closed-loop system with proposed controller\label{system-structure}.}
\end{figure*}

Using the representation of the actuator state \eqref{act}--\eqref{error-control}, the original system \eqref{u-main}--\eqref{initial-u} is equivalent to the following PDE-PDE cascade system
\begin{align} 
\label{cascade-u-main}
u_{t}(x,t)&=u_{xx}(x,t)+\lambda u(x,t)\\
\label{cascade-u-bnd0}
u(0,t)&=0\\
\label{cascade-u-bnd1}
u(1,t)&=\hat{v}(0,t)+\delta(t)^T\bm{\tilde{v}}(0,t)\\
u(x,0)&=u_0(x)\\
\tau_0\hat{v}_t(x,t)&=\hat{v}_{x}(x,t)\label{cascade-hat-v-main}\\
\hat{v}(1,t)&=U(t)\label{cascade-v-control}\\
\hat{v}(x,0)&=\hat{v}_0(x)\label{initial-breve-v0}\\
\Lambda_{\tau}\bm{\tilde{v}}_t(x,t)&=\bm{\tilde{v}}_{x}(x,t)-\Sigma_{\tau}\hat{v}_x(x,t)\label{cascade-ev-main}\\
\bm{\tilde{v}}(1,t)&=\bm{0} \label{cascade-ev-boundary}\\
\bm{\tilde{v}}(x,0)&=\bm{\tilde{v}}_0(x)\label{cascade-initial-tilde-v0}
\end{align}
where $\Sigma_{\tau}=\left(\frac{\tau_1-\tau_0}{\tau_0}, \frac{\tau_2-\tau_0}{\tau_0}, \cdots, \frac{\tau_r-\tau_0}{\tau_0}\right)^T$ and  $\bm{0}$  is the r-by-1 all-zeros vector. The matrices $\Lambda_{\tau}=diag(\tau_1, \tau_2, \cdots, \tau_r)$,  $\delta(t)=(\delta_1(t),\cdots,\delta_j(t),\cdots,$ $\delta_r(t))^T$ $\in\mathbb{R}^r~(j=1,2,\cdots,r)$ are such  that, if $\tau(t)=\tau_i$,
\begin{align}\label{delta}
\delta_j(t)=\left\{
\begin{aligned}
&1,&if~j=i\\
&0,&else.
\end{aligned}\
\right.
\end{align}
Note  that $\delta(t)$ and the Markov process $\tau(t)$  have the same transition probabilities. The structure of the closed-loop system depicted in Fig.\ref{system-structure}.

\subsection{Backstepping transformation and controller design}
We construct a nominal delay-compensated  control law for system \eqref{u-main}--\eqref{initial-u} subject to a nominal delay $\tau_0$ using the backstepping transformations \cite{Krstic2009control}
\begin{align}
w(x,t)=&u(x,t)-\int_{0}^{x}p(x,y)u(y,t)dy\label{transformation-uw}\\
\hat{z}(x,t)=&\hat{v}(x,t)-\int_{0}^{1}\gamma(x,y)u(y,t)dy\nonumber\\
&-\tau_0\int_{0}^{x}k(x,y)\hat{v}(y,t)dy\label{transformation-zv}
\end{align}
together with  the inverse transformations 
\begin{align}
u(x,t)=&w(x,t)+\int_{0}^{x}q(x,y)w(y,t)dy\label{inverse-transformation-wu}\\
\hat{v}(x,t)=&\hat{z}(x,t)+\int_{0}^{1}\eta(x,y)w(y,t)dy\nonumber\\
&+\tau_0\int_{0}^{x}l(x,y)\hat{z}(y,t)dy.\label{inverse-transformation-vz}
\end{align}
Defining  $h(x,t)=\hat{v}_x(x,t)$ and using  \eqref{transformation-uw}--\eqref{inverse-transformation-vz}, the cascade system \eqref{cascade-u-main}--\eqref{cascade-initial-tilde-v0} can be transformed into the following target system
\begin{align} 
\label{target-w-main}
w_{t}(x,t)&=w_{xx}(x,t)\\
\label{target-w-bnd0}
w(0,t)&=0\\
\label{target-w-bnd1}
w(1,t)&=\hat{z}(0,t)+\delta(t)^T\bm{\tilde{v}}(0,t)\\
w(x,0)&=w_0(x)\label{target-initial-w}\\
\tau_0\hat{z}_t(x,t)&=\hat{z}_{x}(x,t)+\tau_0k(x,0)\delta(t)^T\bm{\tilde{v}}(0,t)\label{target-z-main}\\
\hat{z}(1,t)&=0\label{target-z-boundary}\\
\hat{z}(x,0)&=\hat{z}_0(x)\label{target-initial-z}\\
\Lambda_{\tau}\bm{\tilde{v}}_t(x,t)&=\bm{\tilde{v}}_{x}(x,t)-\Sigma_{\tau}h(x,t)\label{target-ev-main}\\
\bm{\tilde{v}}(1,t)&=\bm{0} \label{target-ev-boundary}\\
\bm{\tilde{v}}(x,0)&=\bm{\tilde{v}}_0(x)\label{target-initial-v}\\
\tau_0h_t(x,t)&=h_x(x,t)\label{main-h}\\
h(1,t)&=\hat{z}_x(1,t)+\int_{0}^{1}\eta_x(1,y)w(y,t)dy\nonumber\\
&~~~+\tau_0\int_{0}^{1}l_x(1,y)\hat{z}(y,t)dy\label{target-bound2-h}
\end{align}
where 
\begin{align}\label{h-x-t}
h(x,t)=&\hat{z}_x(x,t)+\int_{0}^{1}\eta_x(x,y)w(y,t)dy+\tau_0l(x,x)\hat{z}(x,t)\nonumber\\
&+\tau_0\int_{0}^{x}l_x(x,y)\hat{z}(y,t)dy.
\end{align}
From  \cite{Krstic2009control} and \cite{Wang2021adaptive}, one can deduce that the gain kernels of transformations \eqref{transformation-uw}, \eqref{transformation-zv} satisfy
\begin{align}
p(x,y)&=-\lambda y\frac{I_1\left(\sqrt{\lambda(x^2-y^2)}\right)}{\sqrt{\lambda(x^2-y^2)}}\label{solution-p}\\
\gamma(x,y)&=\sum_{n=1}^{\infty}e^{\tau_0(\lambda-n^2\pi^2)x}\sin(n\pi y)p_{1n}\label{solution-gamma}\\
k(x,y)&=-\sum_{n=1}^{\infty}(-1)^nn\pi e^{\tau_0(\lambda-n^2\pi^2)(x-y)}p_{1n}\label{solution-k}
\end{align}
where $I_1(\cdot)$ is the first-order modified Bessel function of the first kind, $p_{1n}=2\int_{0}^{1}\sin(n\pi \xi)p(1,\xi)d\xi$, which is the Fourier coefficient of $p(1,y)$, i.e., the function $p(1,y)$ can be represented in the form of Fourier series as $p(1,y)=\sum_{n=1}^{\infty}p_{1n}\sin(n\pi y)$.
Similarly,  the gain kernels of the inverse transformations \eqref{inverse-transformation-wu}, \eqref{inverse-transformation-vz} are given as 
\begin{align}
q(x,y)&=-\lambda y\frac{J_1\left(\sqrt{\lambda(x^2-y^2)}\right)}{\sqrt{\lambda(x^2-y^2)}}\label{solution-q}\\
\eta(x,y)&=\sum_{n=1}^{\infty}e^{-\tau_0n^2\pi^2x}\sin(n\pi y)q_{1n}\label{solution-eta}\\
l(x,y)&=-\sum_{n=1}^{\infty}(-1)^nn\pi e^{-\tau_0n^2\pi^2(x-y)}q_{1n}\label{solution-l}
\end{align}
where $J_1(\cdot)$ is the first-order Bessel function, $q_{1n}=2\int_{0}^{1}\sin(n\pi \xi)q(1,\xi)d\xi$.

Substituting \eqref{cascade-v-control}, \eqref{target-z-boundary}, \eqref{solution-p}--\eqref{solution-k} into \eqref{transformation-zv} results into the following delay-compensated boundary control law
\begin{align*}
U(t)=&2 \sum_{n=1}^{\infty}\int_{0}^{1}\sin(n\pi \xi)\lambda\xi\frac{I_1\left(\sqrt{\lambda(1-\xi^2)}\right)}{\sqrt{\lambda(1-\xi^2)}}d\xi\nonumber \\
&\times\left(-\int_{0}^{1}e^{(\lambda-n^2\pi^2)\tau_0}\int_{0}^{1}  \sin(n\pi y)u(y,t)dy\right.\nonumber\\
&\left.+\tau_0\pi n(-1)^n\int_{0}^{1}   e^{\tau_0\left(\lambda-n^2\pi^2\right)(1-y)}\hat{v}(y,t)dy\right)
\end{align*}

which can be rewritten as 
\begin{align}\label{delay-controller}
U(t)=&2\sum_{n=1}^{\infty}\int_{0}^{1}\sin(n\pi \xi)\lambda\xi\frac{I_1\left(\sqrt{\lambda(1-\xi^2)}\right)}{\sqrt{\lambda(1-\xi^2)}}d\xi\nonumber\\
&\times\left(-\int_{0}^{1}e^{\tau_0(\lambda-n^2\pi^2)}\sin(n\pi y) u(y,t)dy\right.\nonumber\\
&\left.+\tau_0\pi n(-1)^n\int_{t-\tau_0}^{t}e^{\left(\lambda-n^2\pi^2\right)(t-s)} U(s)ds\right).
\end{align}

\begin{Remark}\rm
The controller \eqref{delay-controller} is explicit as the gain kernels $p(x,y)$, $\gamma(x,y)$ and $k(x,y)$ are explicitly given. Note that the transformations \eqref{transformation-uw}, \eqref{transformation-zv} used in this contribution do not differ from the ones introduced in \cite{Krstic2009control}. 
  However, the randomness of delay error disturbance induced  that appears in \eqref{cascade-u-bnd1}, is a disturbance  affecting the nominal  delay $\tau_0$. The
 term $\hat{v}(0,t)$ compensates the effect of  the nominal actuator state $\hat{v}(x,t)$ 
 upon which  \eqref{delay-controller} is constructed. Note that the nominal actuator state is strongly coupled with the stochastic error disturbance estimation $\bm{\tilde{v}}(x,t)$ through  \eqref{cascade-ev-main} (see.  Fig.\ref{system-structure}). 
\end{Remark}

 Based on the full-state feedback law \eqref{delay-controller}, our aim is to prove the stability of the original plant \eqref{u-main}--\eqref{initial-u}, which suffers a random boundary input delay. The well-posedness of \eqref{u-main}--\eqref{initial-u} must be clarified first by examining that of \eqref{target-w-main}--\eqref{h-x-t}. 

\section{Well-posedness Study}\label{6}
The proof of  well-posedness of the target system \eqref{target-w-main}--\eqref{h-x-t} rests on $C_0$-semigroup theory  \cite{zhang2022static,han2016exponential}. Note that if the target system is well-posed, the existence and uniqueness of solutions to the original closed-loop system can be directly established from the invertibility of \eqref{transformation-uw}--\eqref{transformation-zv}. Introducing the following change of variable
\begin{align}\label{trans-m-w}
\hspace{-.25cm}m(x,t)=w(x,t)-x\hat{z}(0,t)-x\delta(t)^T\bm{\tilde{v}}(0,t)
\end{align}
and  considering that the sample realization of Markov process $\delta(t)$ is not only right-continuous but has finite state values, one can deduce the existence of  a sequence $\{t_k, k=1,2,\cdots\}$, such that $t_0=0$, $\lim\limits_{k\rightarrow+\infty}t_k=\infty$, and as $t_k \ \leq t < t_{k+1}$, $\delta(t)=\delta(t_{k})$.	

	We are in a position to prove  iteratively that system \eqref{target-w-main}--\eqref{h-x-t} is well-posed  between successive jump times following \cite{zhang2017stochastic}. 
 Let $T\in \mathbb{R}^+$ be an arbitrary positive constant. Now  considering the cascade system \eqref{target-w-main}--\eqref{h-x-t} on the time interval $t\in[t_0, t_1\wedge T]$, we deduce that, $\delta(t)=\delta(t_0)$.  Assuming that $\delta(t_0)=e_i$,
 system \eqref{target-w-main}--\eqref{h-x-t} can be converted into the following one
\begin{align} 
\label{stable-m-main}
m_{t}(x,t)&=m_{xx}(x,t)-x\hat{z}_t(0,t)-xe_i^T\bm{\tilde{v}}_t(0,t)\\
\label{stable-m-bnd0}
m(0,t)&=0,\quad
m(1,t)=0 \\
m(x,0)&=m_0(x)\label{stable-m-initial}\\
\tau_0\hat{z}_t(x,t)&=\hat{z}_{x}(x,t)+\tau_0k(x,0)e_i^T\bm{\tilde{v}}(0,t)\label{stable-z-main}\\
\hat{z}(1,t)&=0\label{stable-z-boundary}\\
\hat{z}(x,0)&=\hat{z}_0(x)\label{stable-z-initial}\\
\Lambda_{\tau}\bm{\tilde{v}}_t(x,t)&=\bm{\tilde{v}}_{x}(x,t)-\Sigma_{\tau}\hat{h}(x,t)\label{stable-ev-main}\\
\bm{\tilde{v}}(1,t)&=\bm{0}\label{stable-ev-boundary}\\
\bm{\tilde{v}}(x,0)&=\bm{\tilde{v}}_0(x)\label{stable-v-initial}\\
\tau_0\hat{h}_t(x,t)&=\hat{h}_x(x,t)\label{main-hat-h}\\
\hat{h}(1,t)&=\hat{z}_x(1,t)+\int_{0}^{1}\eta_x(1,y)(m(y,t)+y\hat{z}(0,t)\nonumber\\
&~~~+ye_i^T\bm{\tilde{v}}(0,t))dy+\tau_0\int_{0}^{1}l_x(1,y)\hat{z}(y,t)dy\label{bound2-hat-h}
\end{align}
where

\begin{align}\label{hat-h}
&\hat{h}(x,t)=\hat{z}_x(x,t)+\int_{0}^{1}\eta_x(x,y)(m(y,t)+y\hat{z}(0,t)+ye_i^T\nonumber\\
&\cdot\bm{\tilde{v}}(0,t))dy+\tau_0l(x,x)\hat{z}(x,t)+\tau_0\int_{0}^{x}l_x(x,y)\hat{z}(y,t)dy
\end{align}
according to \eqref{h-x-t} and \eqref{trans-m-w}.

Furthermore, define the state-space  
\begin{align}
\mathscr{H}=L^2_E(0,1)\times H^2(0,1)\times \mathbb{H}^1(0,1)\times H^1(0,1)
\end{align}
where
\begin{align}
L^2_E(0,1)&=\{f\in L^2(0,1)|f(0)=f(1)=0\}\\
\mathbb{H}^1(0,1)&=\{\bm{h}(x)=(h_1(x),h_2(x),\cdots,h_r(x))^T\\
&~~~~~|\bm{h}(1)=\bm{0},h_j(x)\in H^1(0,1), j=1,2,\cdots,r\}.\nonumber
\end{align}
For any $\Phi_i=(f_i,g_i,\bm{h_i},d_i)\in\mathscr{H},~i=1,2$, we define  the inner product, 
\begin{align}
\langle \Phi_1,\Phi_2\rangle_{\mathscr{H}}=&\int_{0}^{1}(f_1f_2+\alpha_1\tau_0(1+x)(g_1g_2+g'_1g'_2)\nonumber\\
&+\alpha_2(1+x)(\bm{h_1}^T\bm{h_2}+{\bm{h'_1}}^T\bm{h'_2})\\
&+\alpha_3\tau_0(1+x)d_1d_2+\alpha_4\tau_0(1+x)d'_1d'_2)dx\nonumber
\end{align}
where $\alpha_k$, $k=1,2,3,4$, are positive constants.
Let the operator $\mathscr{A}_i$, $i=1,2,\cdots,r$ in $\mathscr{H}$, with    $(\mathscr{H},\Vert \cdot \Vert_{\mathscr{H}})$ being  a Hilbert space such that
\begin{align}\label{operator}
\mathscr{A}_i\left(
\begin{array}{c}
f \\
g \\
\bm{h}\\
d
\end{array}
\right)
=
\left(
\begin{array}{ccc}
f''+\phi_i(x) \\
\frac{1}{\tau_0}g'+k(x,0)e_i^T\bm{h}(0) \\
\Lambda_{\tau}^{-1}\bm{h}'-\Lambda_{\tau}^{-1}\Sigma_{\tau}d\\
\frac{1}{\tau_0}d'
\end{array}
\right)\end{align}
where \begin{align*}
	\phi_i(x)=&-\frac{1}{\tau_0}xg'(0)-xk(0,0)e^T_i\bm{h}(0)\nonumber\\
	&-xe^T_i\Lambda_{\tau}^{-1}(\bm{h}'
	(0)
	-\Sigma_{{\tau}}d(0)).
\end{align*}
The domain of the operator $\mathscr{A}_i$ is defined as 
\begin{align}
&D(\mathscr{A}_i)=\left\{(f,g,\bm{h},d)^T\in\mathscr{H}\bigg|\right.
\\
&\left.\begin{array}{cc}
f\in H^2(0,1)\cap L^2_E(0,1), g,d\in H^1(0,1), \\
\bm{h}\in \mathbb{H}^1(0,1), f(0)=f(1)=g(1)=0,\bm{h}(1)=\bm{0} \\
d(1)=g'(1)+\int_0^1\eta_x(1,y)(f(y)+yg(0)\\~~~+yh_i(0))dy+\tau_0\int_0^1l_x(1,y)g(y)dy
\nonumber
\end{array}
\right\}.
\end{align}
Using  \eqref{operator},  one can rewrite the system \eqref{stable-m-main}--\eqref{hat-h} as an evolution equation in  $\mathscr{H}$ as follows
\begin{align}
\frac{d\Phi(t)}{dt}&= \mathscr{A}_i \Phi(t),~~~t\in[t_0,t_1\wedge T],~ i=1,2,\cdots,r\label{ab-phi}\\
\Phi(0)&=\Phi_0\label{ab-phi0}
\end{align}
where $\Phi(t)=(m(x,t),\hat{z}(x,t),\bm{\tilde{v}}(x,t),\hat{h}(x,t))^T$, and $\Phi_0=(m_0(x),\hat{z}_0(x),\bm{\tilde{v}}_0(x),\hat{h}_0(x))^T$.

Subsequently,  we state the following  lemma.
\begin{Lemma}\label{le-iso}\rm
	For any $i=1,2,\cdots,r$, consider $\mathscr{A}_i$ defined in \eqref{operator}.  Then,  $\mathscr{A}_i^{-1}$ exists and is compact on $\mathscr{H}$. Moreover, 
 $\sigma(\mathscr{A}_i)$, the spectrum of $\mathscr{A}_i$ consists of all isolated eigenvalues of finite multiplicity.
\end{Lemma}
\begin{Proof}\rm
	Defining  $F_i=(\mu_i,\nu_i,\bm{\omega}, \zeta)^T\in\mathscr{H}$, we need to find the solution to $\mathscr{A}_i\Phi=F_i \in  D(\mathscr{A}_i)$, namely,  $(f,g,\bm{h},d)\in D(\mathscr{A}_i)$, such that the following set differential equation admits a solution
	\begin{align}
	f''+\phi_i(x)&=\mu_i(x)\label{f1}\\
	\frac{1}{\tau_0}g'+k(x,0)e_i^T\bm{h}(0)&=\nu_i(x) \\
	\Lambda_{\tau}^{-1}\bm{h}'-\Lambda_{\tau}^{-1}\Sigma_{\tau}d&=\bm{\omega}(x)\\
	\frac{1}{\tau_0}d'&=\zeta(x)\label{d}
	\end{align}
	when the  boundary conditions are given by
	\begin{align}\label{bf1}
	f(0)=&f(1)=0,\quad
	g(1)=0,\quad
	\bm{h}(1)=\bm{0}\\
	d(1)=&-\tau_0k(1,0)e^T_i\bm{h}(0)+\int_{0}^{1}\eta_x(1,y)(f(y)+yg(0)\nonumber\\
	&+ye^T_i\bm{h}(0))dy+\tau_0\int_{0}^{1}l_x(1,y)g(y)dy.\label{d1}
	\end{align}
	Solving \eqref{f1}--\eqref{d1},  the following set of functions are computed
	\begin{align}
	f(x)=&\int_{0}^{x}\int_{0}^{y}\mu_i(s)dsdy+\frac{\nu_i(0)}{6}x^3+\frac{e_i^T\bm{\omega}(0)}{6}x^3\\
	&-x\int_{0}^{1}\int_{0}^{y}\mu_i(s)dsdy-\frac{\nu_i(0)}{6}x-\frac{e_i^T\bm{\omega}(0)}{6}x\nonumber\\
	g(x)=&-\tau_0\int_{x}^{1}\nu_i(s)ds+\tau_0\int_{x}^{1}k(s,0)dse^T_i\bm{h}(0)\\
	\bm{h}(x)=&-\Lambda_{\tau}\int_{x}^{1}\bm{\omega}(s)ds+\tau_0\Sigma_{\tau}\int_{x}^{1}\int_{s}^{1}\zeta(\xi)d\xi ds\nonumber\\
	&-(1-x)\Sigma_{{\tau}}\left(\int_{0}^{1}\eta_x(1,y)f(y)dy-{\tau}^2_0\right.\nonumber\\
	&\cdot\int_{0}^{1}\int_{y}^{1}l_x(1,y)\nu_i(s)dsdy-\tau_0\int_{0}^{1}\eta_x(1,y)dy\nonumber\\
	&\left.\int_{0}^{1}\nu_i(s)ds\right)-(1-x)\Sigma_{{\tau}}\left(\tau_0k(1,0)\int_{0}^{1}y\right.\nonumber\\
	&\cdot\eta_x(1,y)dy+{\tau}^2_0\int_{0}^{1}\int_{y}^{1}l_x(1,y)k(s,0)dsdy\nonumber\\
	&\left.+\tau_0\int_{0}^{1}\eta_x(1,y)dy\int_{0}^{1}k(s,0)ds\right)e_i^T\bm{h}(0)\\
	d(x)=&-\tau_0\int_{x}^{1}\zeta(s)ds-\tau_0\left(\int_{0}^{1}\eta_x(1,y)ydy\right.\\
	&\left.+\tau_0\int_{0}^{1}\int_{y}^{1}l_x(1,y)\nu_i(s)dsdy\right)+\left(-\tau_0k(1,0)\right.\nonumber\\
	&\left.\cdot\int_{0}^{1}\eta_x(1,y)ydy+\tau_0\int_{0}^{1}\eta_x(1,y)dy\int_{0}^{1}k(s,0)\right.\nonumber\\
	&\left.ds+{\tau}^2_0\int_{0}^{1}\int_{y}^{1}l_x(1,y)k(s,0)dsdy\right)e_i^T\bm{h}(0)\nonumber\\
	&+\int_{0}^{1}\eta_x(1,y)\left(\int_{0}^{y}\int_{0}^{s}\mu_i(\xi)d\xi ds+\frac{y^3-y}{6}\right.\nonumber\\
	&\left.\cdot(\mu(0)+e^T_i\bm{\omega}(0))-\int_{0}^{1}\int_{0}^{s}y\mu_i(\xi)d\xi ds\right)dy\nonumber
	\end{align}
	where
	\begin{align*}\label{h0}
	\bm{h}(0)=&\left(I+(1-x)\left(-\tau_0k(1,0)+\int_{0}^{1}\eta_x(1,y)ydy\right.\right.\nonumber\\
	&+\tau_0\int_{0}^{1}\eta_x(1,y)dy\int_{0}^{1}k(s,0)ds\nonumber\\
	&\left.\left.+{\tau}^2_0\int_{0}^{1}\int_{y}^{1}l_x(1,y)k(s,0)dsdy\right)\Sigma_{\tau}e^T_i\right)^{-1}\nonumber\\
	&\cdot\left(-\Lambda_{\tau}\int_{x}^{1}\bm{\omega}(s)ds+\tau_0\Sigma_{\tau}\int_{x}^{1}\int_{s}^{1}\zeta(\xi)d\xi ds\right.\nonumber\\
	&-(1-x)\left(\int_{0}^{1}\eta_x(1,y)f(y)dy+\tau_0\int_{0}^{1}\eta_x(1,y)dy\right.\nonumber\\
	&\left.\left.\cdot\int_{0}^{1}\nu_i(s)ds+{\tau}^2_0\int_{0}^{1}\int_{y}^{1}l_x(1,y)\nu_i(s)dsdy\right)\Sigma_{\tau}\right).
	\end{align*}
	Note that  $I$ is the $r\times r$ unit matrix.

	Above, we have determined the unique solution  $\Phi=(f, g, \bm{h}, d)\in D(\mathscr{A}_i)$ such that $\mathscr{A}_i\Phi=F_i$ and $\mathscr{A}_i^{-1}$ exists. According to the closed operator theorem, operator $\mathscr{A}_i$ is not only closed and bounded but also has a bounded and closed inverse, that is $0\in\rho(\mathscr{A}_i)$ and $\mathscr{A}_i^{-1}: \mathscr{H}\rightarrow D(\mathscr{A}_i)$ is a bounded linear operator. Since $D(\mathscr{A}_i)\subset \mathscr{H}$, the Sobolev Embedding Theorem asserts that $\mathscr{A}_i^{-1}$ is a compact operator on $\mathscr{H}$. Based on the spectral theory of compact operator, $\sigma(\mathscr{A}_i)$ consists of all isolated eigenvalues of finite multiplicity, which concludes the proof of Lemma  \ref{le-iso}.
\end{Proof}

To ensure that the inner-product $\langle \mathscr{A}_i\Phi, \Phi\rangle$ is bounded, let us first introduce the following two lemmas that provide useful bounds to the derivatives of the gain kernels and their inverses for the proof of Theorem \ref{well-t}.
\begin{Lemma}\label{inequ-p}\rm
	Consider the  function $p(x,y)$ defined in \eqref{solution-p}, for all $x\in[0,1]$, $y\in[0,x)$, there exist positive constants $M_i$, $i=1,2,3,4,5$, such that 
	\begin{align}\Vert p(x,y)\Vert^2_{L^2}&\le M_1,\quad  \Vert p(1,y)\Vert^2_{L^2}\le M_2\\
	\Vert p_y(1,y)\Vert^2_{L^2}&\le M_3, \quad \Vert p_{yy}(1,y)\Vert^2_{L^2}\le M_4\\
	\Vert p_{yyy}(1,y)\Vert^2_{L^2}&\le M_5.
	\end{align}
\end{Lemma}
  
 One can use Parseval's theorem and Fourier series to prove Lemma \ref{inequ-p}. Here, the proof is omitted as the reader can find a similar proof procedure in \cite{Krstic2009control}.
\begin{Lemma}\label{inequ-k0}\rm
	For any $x\in[0,1]$, there exists constants $A_k>0$, $k=1,\cdots,8$, such that the functions $k(x,y)$, $\eta(x,y)$ and $l(x,y)$ defined in \eqref{solution-k}, \eqref{solution-eta} and \eqref{solution-l} respectively, satisfy
	\begin{align} \Vert k(x,0)\Vert^2&\le A_1, \quad \Vert k_x(x,0)\Vert^2\le A_2\\
	|k(0,0)|^2&\le A_3,\quad  \Vert \eta_x(1,x)\Vert^2\le A_4\\
	\Vert l_x(1,x)\Vert^2&\le A_5,\quad  \Vert\eta_{xx}(1,x)\Vert^2\le A_6\\
	|l(x,x)|^2&\le A_7,\quad  \Vert l_{xx}(1,x)\Vert^2\le A_8.\end{align}
\end{Lemma}
 
The proof of Lemma \ref{inequ-k0} can be found in Appendix \ref{proof-k0}.
	
\begin{theorem}\rm\label{well-t}
	For any $\mathscr{A}_i\in \mathscr{H}$, $i=1,2,\cdots,r$, $\mathscr{A}_i$ generates a $C_0$-semigroup on $\mathscr{H}.$
\end{theorem}
\begin{Proof}\rm
First, according to   Poincar\'e inequality, i.e., $\frac{1}{4}\int_{0}^{1}m^2(x,t)dx\le\int_{0}^{1}m^2_x(x,t)dx$, 
 one can get
 \begin{align}
  \frac{1}{4}\int_{0}^{1}f^2dx\le\int_{0}^{1}(f')^2dx.
 \end{align}
Then, according to  \eqref{bf1}, one can state the following inequalities
 \begin{align}
 g^2(0)\le&\int_0^1(g'(x))^2dx\\
 h^2_i(0)\le&\int_0^1(h'_i(x))^2dx,~i=1,2,\cdots,r.
\end{align}
Using the estimates above,  for any real $\Phi=(f,g,\bm{h},d)^T\in D(\mathscr{A}_i)$, $i=1,2,\cdots,r$, applying the Cauchy-Schwarz inequality,  Lemma \ref{inequ-k0}, and Lemma \ref{inequs-inner} (see Appendix \ref{A-Inner}), 
the following estimate holds:
	\begin{align}
	&\langle \mathscr{A}_i\Phi, \Phi\rangle
 \le-\left(\frac{1}{4}-\frac{\gamma_1}{2}-5A_4\left(\alpha_3+\alpha_2\frac{r(\max\delta_{\tau})^2}{\underline{\tau}^2\tau^2_0}\right)\right.\nonumber\\
	&\left.-8\alpha_4A_6\right)\int_{0}^{1}f^2dx-\left(-\alpha_1-\frac{5A_4\alpha_2r(\max\delta)^2}{3\underline{\tau}\tau_0^2}-\frac{5A_4\alpha_3}{3}\right.\nonumber\\
 &\left.-A_2\tau_0\gamma_2\alpha_1\right)\int_{0}^{1}(g')^2dx-\left(\frac{\alpha_1}{2}-\alpha_1A_1\tau_0\gamma_2-5A_5\right.\nonumber\\
 &\left.\cdot\left(\alpha_3\tau_0^2+\frac{\alpha_2 r(\max\delta_{\tau})^2}{\underline{\tau}}\right)-8\tau^2_0A_8\alpha_4\right)\int_{0}^{1}g^2dx\nonumber\\
	&-\alpha_2\left(\frac{1}{2\overline{\tau}}-\frac{(\max\delta_{\tau})^2}{2\gamma_3\tau^2_0\underline{\tau}^2}\right)\sum_{j=1}^{r}\int_{0}^{1}h_j^2dx-\left(-\frac{\alpha_2}{\overline{\tau}}-\frac{2A_1}{3\gamma_1}\right.\nonumber\\
	&-\frac{\alpha_2(\max\delta_{\tau})^2}{2\gamma_4\tau^2_0\underline{\tau}^2}-\frac{7\tau_0\alpha_1}{6\gamma_2}-\alpha_1\tau_0^2A_1-\left(5\tau_0^2A_1+\frac{A_4}{3}\right)\nonumber\\
 &\cdot\left(\alpha_3+\frac{r(\max\delta_{\tau})^2\alpha_2}{\underline{\tau}^2\tau_0^2}\right)-8\alpha_4\left(\tau_0^2A_2+\tau_0^4A_1A_7\right.\nonumber\\
	&\left.\left.+\frac{A_6}{3}\right)\right)\sum_{j=1}^{r}\int_{0}^{1}(h'_j)^2dx-\left(\frac{\alpha_3}{2}-2\gamma_3\alpha_2\right)\int_{0}^{1}d^2dx\nonumber\\
 &-\left(\frac{\alpha_4}{2}-2\gamma_4\alpha_2\right)\int_{0}^{1}(d')^2dx-\left(\frac{\alpha_1}{2}-\frac{2}{3\gamma_1\tau^2_0}\right)\nonumber\\
 &\cdot(g'(0))^2-\left(\frac{\alpha_2}{2\overline{\tau}}-\frac{2}{3\gamma_1\underline{\tau}^2}-\frac{8\tau^4_0A_1\alpha_4}{\underline{\tau}^2}\right)(h'_i(0))^2
	\nonumber\\
	&-\left(\frac{\alpha_3}{2}-\frac{8(\max\delta_{\tau})^2\tau^2_0A_1\alpha_4}{\underline{\tau}^2}-\frac{2(\max\delta_{\tau})^2}{3\gamma_1\underline{\tau}^2\tau_0^2}\right)d^2(0)\nonumber\\
 &-\frac{\alpha_4}{2}(d'(0))^2-\frac{\alpha_2}{2\overline{\tau}}\sum_{j=1,j\ne i}^{r}\left(h_j^2(0)+(h'_j(0))^2\right)\nonumber\\
	&\le M\langle \Phi,\Phi\rangle_{\mathscr{H}}\label{dissip}
	\end{align}
	where $\max\delta_{\tau}=\max|\tau_0-\tau_i|$, $i=1,2,\cdots,r$, 
	$\gamma_k,~k=1,2,3,4$, are positive constants,
\begin{align}
\alpha_1>&\frac{4}{3\gamma_1\tau_0^2}\\
\alpha_2>&\frac{4\overline{\tau}}{3\gamma_1\underline{\tau}^2}+\frac{16\tau_0^4\overline{\tau}A_1\alpha_4}{\underline{\tau}^2}\\
\alpha_3>&\frac{16(\max\delta_{\tau})^2\tau^2_0A_1\alpha_4}{\underline{\tau}^2}+\frac{4(\max\delta_{\tau})^2}{3\gamma_1\underline{\tau}^2\tau_0^2}\\
\alpha_4>&0.
\end{align}
  Since $\max\delta_{\tau}$ is bounded, the constant in \eqref{dissip} holds can be selected as follows \begin{align}
     M=&\max\left\{\left|\frac{1}{4}-\frac{\gamma_1}{2}-5A_4\left(\alpha_3+\alpha_2\frac{r(\max\delta_{\tau})^2}{\underline{\tau}^2\tau^2_0}\right)\right.\right.\nonumber\\
	&\left.-8\alpha_4A_6\right|,\left|-\alpha_1-\frac{5A_4\alpha_2r(\max\delta)^2}{3\underline{\tau}\tau_0^2}-\frac{5A_4\alpha_3}{3}\right.\nonumber\\
 &\left.-A_2\tau_0\gamma_2\alpha_1\right|,\left|\frac{\alpha_1}{2}-\alpha_1A_1\tau_0\gamma_2-\left(\frac{\alpha_2 r(\max\delta_{\tau})^2}{\underline{\tau}}\right.\right.\nonumber\\
 &\left.\left.+\alpha_3\tau_0^2\right)\cdot5A_5-8\tau^2_0A_8\alpha_4\right|,\alpha_2\left|\frac{1}{2\overline{\tau}}-\frac{(\max\delta_{\tau})^2}{2\gamma_3\tau^2_0\underline{\tau}^2}\right|,\nonumber\\
 &\left|-\frac{\alpha_2}{\overline{\tau}}-\frac{2A_1}{3\gamma_1}-\frac{\alpha_2(\max\delta_{\tau})^2}{2\gamma_4\tau^2_0\underline{\tau}^2}-\frac{7\tau_0\alpha_1}{6\gamma_2}-\alpha_1\tau_0^2A_1\right.\nonumber\\
 &-\left(\alpha_3+\frac{r(\max\delta_{\tau})^2\alpha_2}{\underline{\tau}^2\tau_0^2}\right)\left(5\tau_0^2A_1+\frac{A_4}{3}\right)\nonumber\\
 &\left.-8\alpha_4\left(\tau_0^2A_2+\frac{A_6}{3}+\tau_0^4A_1A_7\right)\right|,\left|\frac{\alpha_3}{2}-2\gamma_3\alpha_2\right|,\nonumber\\
 &\left.\left|\frac{\alpha_4}{2}-2\gamma_4\alpha_2\right|\right\}.\nonumber
 \end{align}
	Knowing that $\langle \mathscr{A}_i\Phi, \Phi\rangle\le M\langle \Phi,\Phi\rangle_{\mathscr{H}},$  we deduce that $\mathscr{A}_i-MI$ is a dissipative operator. Based on Lemma \ref{le-iso}, we can get that $\mathscr{A}_i-MI$ satisfies the conditions of Lumer-Phillips theorem (see.  page 14 in Reference \cite{pazy2012semigroups})  and therefore, $\mathscr{A}_i$ generates a $C_0$-semigroup $e^{\mathscr{A}_it}$ on $\mathscr{H}$. The well-posedness of the system \eqref{ab-phi}--\eqref{ab-phi0} gets proven. Because  $T$ is arbitrarily defined, the above process can be repeated iteratively to prove that system \eqref{ab-phi}--\eqref{ab-phi0} is well-posed for any $t>0$.
\end{Proof}

Based on the invertibility of the backstepping transformations \eqref{transformation-uw}--\eqref{transformation-zv}, combining the Lemma \ref{le-iso} and  Theorem \ref{well-t} which state the well-posedness of the target system, we can infer that the cascade system \eqref{cascade-u-main}--\eqref{cascade-initial-tilde-v0} under the controller \eqref{delay-controller} is well-posed as well, so is the original system \eqref{u-main}--\eqref{initial-u}.

\section{Stability Analysis}\label{3}
In this section, we state the stability of cascade system \eqref{cascade-u-main}--\eqref{cascade-initial-tilde-v0} and that of the target system \eqref{target-w-main}--\eqref{h-x-t}.

\begin{theorem}\label{stability-cascade}\rm	
	Consider the closed-loop system \eqref{cascade-u-main}--\eqref{cascade-initial-tilde-v0} under controller \eqref{delay-controller}. For any initial conditions  $(u_0(x), \hat{v}_0(x), \bm{\tilde{v}}_0(x))\in L^2(0,1)\times H^2(0,1)\times\mathbb{H}^1(0,1)$, and compatible conditions satisfying $u_0(0)=0$, $\bm{\tilde{v}}_0(1)=\bm{0}$, there exists a constant $\varepsilon>0$ whose value depends on the upper bound of the transition rate  $q^\star$ and the reactivity coefficient $\lambda$, such that, if 
 \begin{align}\label{consition1-stable}
    \mathbb{E}_{[0,\tau(0)]}\left(|\tau_0-\tau(t)|^2\right)\le\varepsilon (\lambda, q^\star)
 \end{align}
there exist positive constants $\alpha$ and $\beta$, such that
	\begin{align}
	\mathbb{E}_{[0,(V_1(0),\tau(0))]}(V_1(t))\le\alpha V_1(0)e^{-\beta t}
	\end{align}
 with
	\begin{align}
	V_1(t)=\Vert u\Vert^2_{L^2}+\Vert \hat{v}\Vert^2_{H^1}+\Vert \bm{\tilde{v}}\Vert^2_{H^1}+\Vert \hat{v}_x\Vert^2_{H^1}
	\end{align}
	where $\mathbb{E}_{[0,(V_1(0),\tau(0))]}(V_1(t))$ is the condition expectation of $V_1(t)$ at instant $t$ under the conditions $V_1(t)=V_1(0)$ and $\tau(t)=\tau(0)$ at $t=0$. \end{theorem}

Theorem \ref{stability-cascade} establishes mean-square exponential stability of system \eqref{cascade-u-main}--\eqref{cascade-initial-tilde-v0}, which implies the stability of the original system \eqref{u-main}--\eqref{initial-u} under the controller \eqref{delay-controller}. To prove Theorem \ref{stability-cascade}, we will first state the norm-equivalence between system \eqref{cascade-u-main}--\eqref{cascade-initial-tilde-v0} and system  \eqref{target-w-main}--\eqref{h-x-t} before performing a Lyapunov analysis of the target system \eqref{target-w-main}--\eqref{h-x-t}.

Let us first state the following two lemmas that show the boundedness of transformations \eqref{transformation-uw}--\eqref{inverse-transformation-vz}. 
\begin{Lemma}\label{boundness-transformation}\rm
	For any given functions $f(x),~g(x)\in L^2(0,1)$, there exist constants $B_i>0$, $i=\{1,\cdots,6\}$, such that
	\begin{align}
	\int_{0}^{1}\left(\int_{0}^{x}p(x,y)f(y)dy\right)^2dx&\le B_1\Vert f\Vert^2_{L^2}\label{inequ-pf}\\
	\int_{0}^{1}\left(\int_{0}^{1}\gamma(x,y)f(y)dy\right)^2dx&\le B_2\Vert f\Vert^2_{L^2}\label{inequ-gammaf}\\
	\int_{0}^{1}\left(\int_{0}^{x}k(x,y)g(y)dy\right)^2dx&\le B_3 \Vert g\Vert^2_{L^2}\label{inequ-kg}\\
	\int_{0}^{1}\left(\int_{0}^{1}\gamma_x(x,y)f(y)dy\right)^2dx&\le B_4\Vert f\Vert^2_{L^2}\label{inequ-gamma-x-f}\\
	\int_{0}^{1}\left(k(x,x)g(x)\right)^2dx&\le B_5 \Vert g_x\Vert^2_{L^2}\label{inequ-k-x-f}\\
	\int_{0}^{1}\left(\int_{0}^{x}k_x(x,y)g(y)dy\right)^2dx&\le B_6 \Vert g\Vert^2_{L^2}.\label{inequ-kg-x}
	\end{align}
\end{Lemma}

\begin{Proof}\rm
	Inequality \eqref{inequ-pf} can be established using  Cauchy-Schwarz inequality and Lemma \ref{inequ-p}.
	
	Applying the Cauchy-Schwarz inequality, Parseval's theorem and Lemma \ref{inequ-p}, \eqref{inequ-gammaf} can be obtained as follows
	\begin{align*}
	&\int_{0}^{1}\left(\int_{0}^{1}\gamma(x,y)f(y)dy\right)^2dx\\
	\le&\int_{0}^{1}\int_{0}^{1}\left(\sum_{n=1}^{\infty}e^{\tau_0(\lambda-n^2\pi^2)x}\sin(n\pi y)p_{1n}\right)^2dydx\Vert f\Vert^2_{L^2}\\
	\le&\sum_{n=1}^{\infty}\frac{e^{2\lambda\tau_0}}{2\tau_0(\lambda-n^2\pi^2)}\sum_{n=1}^{\infty}p_{1n}^2 \Vert f\Vert^2_{L^2}\le B_2  \Vert f\Vert^2_{L^2}.
	\end{align*}
	Again, using the Cauchy-Schwarz inequality and Lemma \ref{inequ-p}, we get 
	\begin{align*}
	&\int_{0}^{1}\left(\int_{0}^{x}k(x,y)g(y)dy\right)^2dx\\
	\le&\int_{0}^{1}\int_{0}^{x}\left(\sum_{n=1}^{\infty}(-1)^nn\pi e^{\tau_0(\lambda-n^2\pi^2)(x-y)}p_{1n}\right)^2dydx\Vert g\Vert^2_{L^2}\\
	\le&\int_{0}^{1}\int_{0}^{x}\sum_{n=1}^{\infty} e^{2\tau_0(\lambda-n^2\pi^2)(x-y)}\sum_{n=1}^{\infty}(n\pi)^2p_{1n}^2dydx\Vert g\Vert^2_{L^2}\\
	\le&\sum_{n=1}^{\infty}\frac{2e^{2\lambda\tau_0}}{\tau_0(\lambda-n^2\pi^2)}\Vert p_{\xi}(1,\xi)\Vert^2_{L^2}\Vert g\Vert^2_{L^2}\le B_3\Vert g\Vert^2_{L^2} 
	\end{align*}
as in \eqref{inequ-kg}.

	By a similar approach, we can prove the inequalities \eqref{inequ-gamma-x-f}--\eqref{inequ-kg-x},  which concludes the proof.
\end{Proof}
\begin{Lemma}\label{inverse-boundness-transformation}\rm
	For any given functions $f(x),~g(x)\in L^2(0,1)$, there exist constants $C_i>0$, $i=\{1,\cdots,6\}$, such that
	\begin{align}
	\int_{0}^{1}\left(\int_{0}^{x}q(x,y)f(y)dy\right)^2dx&\le C_1\Vert f\Vert^2_{L^2}\label{inequ-qf}\\
	\int_{0}^{1}\left(\int_{0}^{1}\eta(x,y)f(y)dy\right)^2dx&\le C_2\Vert f\Vert^2_{L^2}\label{inequ-etaf}\\
	\int_{0}^{1}\left(\int_{0}^{x}l(x,y)g(y)dy\right)^2dx&\le C_3 \Vert g\Vert^2_{L^2}\label{inequ-lg}\\
	\int_{0}^{1}\left(\int_{0}^{1}\eta_x(x,y)f(y)dy\right)^2dx&\le C_4\Vert f\Vert^2_{L^2}\label{inequ-eta-x-f}\\
	\int_{0}^{1}\left(l(x,x)g(x)\right)^2dx&\le C_5 \Vert g_x\Vert^2_{L^2}\label{inequ-l-x-f}\\
	\int_{0}^{1}\left(\int_{0}^{x}l_x(x,y)g(y)dy\right)^2dx&\le C_6 \Vert g\Vert^2_{L^2}.\label{inequ-lg-x}
	\end{align}
\end{Lemma}

This lemma can be  proved  proceeding with the proof methodology of   Lemma \ref{boundness-transformation}.

\begin{Proposition}\label{prop2}\rm
	There exist positive constants $\beta_1$ and $\beta_2$, such that
	\begin{align}
	&\beta_1(\Vert u\Vert^2_{L^2}+\Vert \hat{v}\Vert^2_{H^1}+\Vert \bm{\tilde{v}}\Vert^2_{H^1}+\Vert \hat{v}_x\Vert^2_{H^1})\nonumber\\
	\le&(\Vert w\Vert^2_{L^2}+\Vert \hat{z}\Vert^2_{H^1}+\Vert \bm{\tilde{v}}\Vert^2_{H^1}+\Vert h\Vert^2_{H^1})\nonumber\\
	\le&\beta_2(\Vert u\Vert^2_{L^2}+\Vert \hat{v}\Vert^2_{H^1}+\Vert \bm{\tilde{v}}\Vert^2_{H^1}+\Vert \hat{v}_x\Vert^2_{H^1}).\label{inequ-target-cascade}
	\end{align}
\end{Proposition}
\begin{Proof}\rm
	Based on the transformations \eqref{transformation-uw}, \eqref{transformation-zv} and the fact that $\hat{v}_x(x,t)=h(x,t)$, one can use   Lemma \ref{boundness-transformation} and Cauchy-Schwarz inequality to  infer that
	\begin{align}
	&\Vert w\Vert^2_{L^2}+\Vert \hat{z}\Vert^2_{H^1}+\Vert \bm{\tilde{v}}\Vert^2_{H^1}+\Vert h\Vert^2_{H^1}\nonumber\\
	\le&(2+2B_1+3B_2+4B_4)\Vert u\Vert^2_{L^2}+4(1+\tau^2_0B_5)\Vert\hat{v}_x\Vert^2_{L^2}\nonumber\\
	&+(3+3B_3\tau^2_0+4\tau^2_0B_6)\Vert\hat{v}\Vert^2_{L^2}+\Vert \bm{\tilde{v}}\Vert^2_{H^1}+\Vert \hat{v}_x\Vert^2_{H^1}\nonumber\\
	\le&\beta_2(\Vert u\Vert^2_{L^2}+\Vert \hat{v}\Vert^2_{H^1}+\Vert \bm{\tilde{v}}\Vert^2_{H^1}+\Vert \hat{v}_x\Vert^2_{H^1})
	\end{align}
	where $\beta_2=\max\{2+2B_1+3B_2+4B_4,4(1+\tau^2_0B_5),3+3B_3\tau^2_0+4\tau^2_0B_6\}$. The right of the inequality \eqref{inequ-target-cascade} gets proven. Similarly, using  Lemma \ref{inverse-boundness-transformation}, one can prove the left inequality of \eqref{inequ-target-cascade}.
\end{Proof}

We are in the position to prove the exponential stability of target system \eqref{target-w-main}--\eqref{h-x-t} by first introducing the following lemma to obtain some useful bounds of the gain kernels and their derivatives.
	\begin{Lemma}\label{inequ-k1-gamma1}\rm
	The following bounds can be obtained for the functions $\gamma(x,y)$ and $k(x,y)$  defined in \eqref{solution-gamma} and \eqref{solution-k}
 for all $0\leq y\leq x\leq 1$, 
	\begin{align}
	\Vert \gamma(0,y)\Vert^2,\Vert \gamma(1,y)\Vert^2&\le D_1\\
 \Vert \gamma_y(0,y)\Vert^2,\Vert \gamma_y(1,y)\Vert^2&\le D_2\\
	\Vert \gamma_{yy}(0,y)\Vert^2,\Vert \gamma_{yy}(1,y)\Vert^2&\le D_3\\
 \Vert k(1,y)\Vert^2\le D_4,~~
	\Vert k_{y}(1,y)\Vert^2&\le D_5\\
	\Vert\gamma_{yyyy}(1,y)\Vert^2\le D_6,~~|\gamma_{yyy}(1,1)|^2&\le D_7\\
 	\Vert k_{yy}(1,y)\Vert^2&\le D_8
	\end{align}
	where  $D_i>0$, $i=\{1,\cdots,8\}$ are constants.
\end{Lemma}

Following the steps of the proof of Lemma \ref{inequ-k0}, one can prove Lemma \ref{inequ-k1-gamma1}.
\begin{Proposition}\label{stbility-tsrget-m}\rm
	Consider the target system \eqref{target-w-main}--\eqref{h-x-t}. For any initial conditions  $(w_0(x), \hat{z}_0(x), \bm{\tilde{v}}_0(x))\in L^2(0,1)\times H^2(0,1)\times\mathbb{H}^1(0,1)$, and compatible conditions satisfying $w_0(0)=0$, $\hat{z}_0(1)=0$, $\bm{\tilde{v}}_0(1)=\bm{0}$, there exist constants $\hat{\alpha}_1, \hat{\beta}_1>0$, and $\varepsilon>0$ whose value depends on the upper bound of the transition rate  $q^\star$ and the reactivity coefficient $\lambda$, satisfying
 \begin{align}\label{consition2-stable}
    \mathbb{E}_{[0,\tau(0)]}\left(|\tau_0-\tau(t)|^2\right)\le\varepsilon (\lambda, q^\star)
 \end{align} 
 such that
	\begin{align}
	\mathbb{E}_{[0,(V_2(0),\tau(0))]}(V_2(t))\le\hat{\alpha}_1 V_2(0)e^{-\hat{\beta}_1 t}
	\end{align}
	where 
	\begin{align}
	&V_2(t, \delta(t))
	=\int_{0}^{1}w^2(x,t)dx+a\tau_0\int_{0}^{1}\left(\hat{z}^2(x,t)\right.\nonumber\\
	&\left.+\hat{z}^2_x(x,t)\right)dx+b\int_{0}^{1}\left(\left(\delta^T(t)\Lambda_{\tau}\bm{\tilde{v}}_x(x,t)\right)^2\right.\\
	&\left.+\left(\delta^T(t)\Lambda_{\tau}\bm{\tilde{v}}(x,t)\right)^2\right)dx+c\int_{0}^{1}\left(h^2(x,t)+h^2_x(x,t)\right)dx\nonumber
	\end{align}
	with $a,b,c$ being positive constants.
\end{Proposition}
\begin{Proof}\rm
	Define the Lyapunov candidate
	\begin{align}
	&V_3(t, \delta(t))
	=\int_{0}^{1}m^2(x,t)dx+\int_{0}^{1}(1+x)\left[a\tau_0\left(\hat{z}^2(x,t)\right.\right.\nonumber\\
	&\left.+\hat{z}^2_x(x,t)\right)+b\left(\left(\delta^T(t)\Lambda_{\tau}\bm{\tilde{v}}(x,t)\right)^2+\left(\delta^T(t)\Lambda_{\tau}\bm{\tilde{v}}_x(x,t)\right)^2\right)\nonumber\\
	&\left.+c\tau_0\left(\hat{h}^2(x,t)+\hat{h}^2_x(x,t)\right)\right]dx\label{V-2}
	\end{align}
	where $a,b,c$ are positive constants. 
	 Clearly, there exist positive constants $\beta_{3}$ and $\beta_{4}$, such that 
  \begin{align}\label{equ-V-2-3}
      \beta_{3} V_2\le V_3\le \beta_{4} V_2.
  \end{align}

	Knowing that $\delta(t)$ is a Markov process,  $V_3(t,\delta(t))$ is not differentiable with respect to time $t,$ we  introduce  the infinitesimal generator $\mathcal{L}$  of $V_3(t,\delta(t))$ \cite{kong2022probabilistic3}: 
 
	\begin{align}\label{operator-L-V3}	\mathcal{L}V_3(t,\delta(t))=&\lim\limits_{\Delta t\rightarrow0^+}\sup\frac{1}{\Delta t}\left(\mathbb{E}\{V_3((t+\Delta t),\delta(t+\Delta t))\}\right.\nonumber\\
	&\left.-V_3(t,\delta(t))\right).
	\end{align}
	For  simplicity, we denote $V_3(t,\delta(t))$ as $V_3(t)$, $V_3(t,\delta(t)=e_i)$ as $V_3(t,i)$. In addition, we define, $\mathcal{L}_i,$ the infinitesimal generator of Markov process $(m, {\hat{z}},\bm{\tilde{v}}, \hat{h},\delta(t)):$ 
	\begin{align*}
\mathcal{L}_iV_3(t)=\frac{dV_3(t,i)}{dt}+\sum_{j=1}^{r}q_{ij}(t)(V_3(t,j)-V_3(t,i)).
	\end{align*}
	Knowing (see. \cite{kong2022probabilistic3})	\begin{align}\label{operator-L-simple}	\mathcal{L}V_3(t)=&\sum_{j=1}^{r}P_{ij}(0,t)\mathcal{L}_jV_3(t).
	\end{align}
	Based on \eqref{operator-L-simple}, we first focus on $\mathcal{L}_iV_3$. Utilizing integration by parts, Cauchy-Schwarz inequality, Young's inequality, Lemma \ref{inequs-V}, Lemma \ref{bound-V3-ij} (see. Appendix \ref{A-h}) and the Poincar\'{e} inequality, i.e., 
 
 \begin{align}
     \frac{1}{4}\int_{0}^{1}m^2(x,t)dx\le\int_{0}^{1}m^2_x(x,t)dx
 \end{align}
 one gets
 
	\begin{align}\label{inequ-part1-v3}
	&\mathcal{L}_iV_3(t)\le-\left(\frac{1}{2}-\gamma_5-5A_4\left(2c+\frac{2\overline{\tau}(\tau_0-\tau_i)^2b}{\tau^2_0}\right)\right.\nonumber\\
	&\left.-16A_6c\right)\int_{0}^{1}m^2(x,t)dx-\left(a-a\gamma_6-5\tau_0^2A_5\right.\nonumber\\
	&\left.\cdot\left(2c+\frac{2\overline{\tau}(\tau_0-\tau_i)^2b}{\tau^2_0}\right)-16A_8\tau_0^2c\right)\int_{0}^{1}\hat{z}^2(x,t)dx\nonumber\\
	&-a\left(1-\gamma_7\right)\int_{0}^{1}\hat{z}^2_x(x,t)dx-b\left(\underline{\tau}-4\gamma_8\right)\int_{0}^{1}\tilde{v}^2_i(x,t)dx\nonumber\\
	&-b\left(\underline{\tau}-4\gamma_{9}\right)\int_{0}^{1}(\partial_x\tilde{v}_i(x,t))^2dx-\left(c-b\frac{\overline{\tau}^2(\tau_0-\tau_i)^2}{\gamma_8{\tau}^2_0}\right.\nonumber\\
 &\left.-\frac{4\overline{\tau}^2b}{\tau_0^2}\sum_{j=1}^{r}|q_{ij}(t)|\left(|\tau_0-\tau_i|^2+|\tau_0-\tau_j|^2\right)\right)\nonumber\\
	&\cdot\int_{0}^{1}\hat{h}^2(x,t)dx-\left(c-\frac{\overline{\tau}^2(\tau_0-\tau_i)^2b}{\gamma_9\tau^2_0}\right)\int_{0}^{1}\hat{h}^2_x(x,t)dx\nonumber\\
 &-\left(a-\frac{4}{3\tau_0^2\gamma_5}-\frac{120A_1\overline{\tau}^4b}{\tau^2_0}\sum_{j=1}^{r}|q_{ij}(t)|\left(|\tau_0-\tau_i|^2\right.\right.\nonumber\\
 &\left.\left.+|\tau_0-\tau_j|^2\right)\right)\hat{z}^2_x(0,t)-\left(a-\frac{16A_6}{3}c-\frac{5A_4}{3}\left(2c\right.\right.\nonumber\\
	&\left.\left.+\frac{2\overline{\tau}(\tau_0-\tau_i)^2b}{\tau^2_0}\right)-14A_1b\overline{\tau}^2\sum_{j=1}^{r}|q_{ij}(t)|\left(|\tau_0-\tau_i|^2\right.\right.\nonumber\\
 &\left.\left.+|\tau_0-\tau_j|^2\right)(1+30D_4\overline{\tau}^2)\right)\hat{z}^2(0,t)-\left(b\underline{\tau}-\frac{4}{3\gamma_5\underline{\tau}^2}\right.\nonumber\\
	&-\frac{16\tau^4_0A_1c}{\underline{\tau}^2}-\frac{120D_2\overline{\tau}^4b}{\tau_0^2}\sum_{j=1}^{r}|q_{ij}(t)|\left(|\tau_0-\tau_i|^2\right.\nonumber\\
 &\left.\left.+|\tau_0-\tau_j|^2\right)\right)(\partial_x\tilde{v}_i(0,t))^2-\left(b\underline{\tau}-a\tau_0^2\left(2A_1+\frac{4A_1}{\gamma_6}\right.\right.\nonumber\\
 &\left.+\frac{4A_2}{\gamma_7}\right)-\left(5\tau^2_0A_1+\frac{5A_4}{3}\right)\left(2c+\frac{2\overline{\tau}(\tau_0-\tau_i)^2b}{\tau^2_0}\right)\nonumber\\
	&-\frac{4A_3}{3\gamma_5}-16\left(\tau^2_0A_2+\frac{A_6}{3}+\tau^4_0A_1A_7\right)c-120b\overline{\tau}^4\nonumber\\
 &\left.\cdot\sum_{j=1}^{r}|q_{ij}(t)|(|\tau_0-\tau_i|^2+|\tau_0-\tau_j|^2)\left(D_7+\lambda^2D_1\right.\right.\nonumber\\
 &\left.\left.+D_2A_1+A_1A_3\right)\right)\tilde{v}_i^2(0,t)-\left(c-\frac{4(\tau_0-\tau_i)^2}{3\gamma_5\tau^2_0\underline{\tau}^2}\right.\nonumber\\
	&\left.-\frac{16\tau^4_0A_1(\tau_0-\tau_i)^2}{\tau^2_0\underline{\tau}^2}c\right)\hat{h}^2(0,t)-c\hat{h}^2_x(0,t)+2b\overline{\tau}^2\nonumber\\
  &\sum_{j=1}^{r}|q_{ij}(t)|\left(|\tau_0-\tau_i|^2+|\tau_0-\tau_j|^2\right)\left(\frac{30\overline{\tau}^2M_{v2_0}}{\beta_3}\right.\nonumber\\
  &\left.+\frac{30\overline{\tau}^2(D_4M_{v1}+M_{v2})+M_{v1}}{\beta_1\beta_3}\right)V_3(t)
	\end{align}
		with
  \begin{align*}
    M_{v1}=&7\max\left\{2D_2,D_3+\lambda^2D_1+2D_1D_4+A_1M_1,\right.\nonumber\\
    &\left.2\tau_0^2D_4^2+D_5\right\}\\
    M_{v2}=&\max\left\{D_6+4\lambda^2D_3+\lambda^4D_1+A_1D_3+\lambda^2A_1D_1\right.\nonumber\\
    &\left.+\frac{2D_1D_5}{\tau_0^2}, 2D_4D_5+\frac{D_8}{\tau_0^2}\right\}\\
    M_{v2_0}=&\max\left\{\frac{D_2}{\tau^2_0}+\frac{2D_2(\tau_0-\tau_i)^2}{\underline{\tau}^2\tau_0^2}, 2\left(D_7+\lambda^2D_1\right.\right.\nonumber\\
    &\left.\left.+D_2A_1+\frac{D_5}{2\tau_0^2}\right)\right\}.
  \end{align*}
	Let 
 \begin{align*}
     0&<\gamma_5<\frac{1}{2}-(10A_4+16A_6)c\\
     0&<\gamma_6,\gamma_7< 1\\
     0&<\gamma_8,\gamma_9<\frac{\underline{\tau}}{4}\\
     a&>\max\left\{\frac{(10A_5+16A_8)\tau_0^2c}{1-\gamma_6},\frac{4}{3\tau_0^2\gamma_5}, \frac{10A_4+16A_6}{3}c\right\}\\
     b&>\frac{1}{\underline{\tau}}\max\left\{\frac{4}{3\gamma_5\underline{\tau}^2}+\frac{16\tau^4_0A_1c}{\underline{\tau}^2},\frac{4A_3}{3\gamma_5}+\frac{a\tau_0^2}{\underline{\tau}}\left(2A_1\right.\right.\nonumber\\
     &~~~\left.+\frac{4A_1}{\gamma_6}+\frac{4A_2}{\gamma_7}\right)+2c\left(5\tau_0^2A_1+8\tau_0^2A_2+\frac{8A_6+5A_4}{3}\right.\nonumber\\
     &\left.\left.~~~+8\tau_0^4A_1A_7\right)\right\}\\
     0&<c<\frac{1}{2(10A_4+16A_6)}
 \end{align*}
and
 \begin{align}
     \theta_1=&\min\left\{\frac{1}{2}-\gamma_5-10A_4c-16A_6c,a-a\gamma_6-10\tau_0^2A_5c\right.\nonumber\\
     &\left.-16A_8\tau_0^2c,a(1-\gamma_7),b(\underline{\tau}-4\gamma_8),b(\underline{\tau}-4\gamma_9),c\right\}\label{theta1}\\
     \theta_2=&\min\left\{a-\frac{4}{3\tau_0^2\gamma_5},a-\frac{16A_6c}{3}-\frac{10A_4c}{3},b\underline{\tau}-\frac{4}{3\gamma_5\underline{\tau}^2}\right.\nonumber\\
	&-\frac{16\tau^4_0A_1c}{\underline{\tau}^2},b\underline{\tau}-\frac{4A_3}{3\gamma_5}-a\tau_0^2\left(2A_1+\frac{4A_1}{\gamma_6}+\frac{4A_2}{\gamma_7}\right)\nonumber\\
 &-2c\left(5\tau_0^2A_1+\frac{5A_4}{3}\right)-16\left(\tau_0^2A_2+\frac{A_6}{3}\right.\nonumber\\
 &\left.\left.+\tau_0^4A_1A_7\right) c,c\right\}\label{theta2}\\
     \rho_1=&\max\left\{\frac{10\overline{\tau}A_4b}{\tau_0^2},\frac{10\overline{\tau}A_5b}{\tau_0^2},\frac{\overline{\tau}^2b}{\gamma_8\tau_0^2},\frac{\overline{\tau}^2b}{\gamma_9\tau_0^2}\right\}\label{rho1}\\
     \rho_2=&\max\left\{\frac{10A_4\overline{\tau}b}{3\tau_0^2},\frac{2\overline{\tau}b}{\tau_0^2}\left(5\tau_0^2A_1+\frac{5A_4}{3}\right),\frac{4}{3\gamma_5\tau_0^2\underline{\tau}^2}\right.\nonumber\\
     &\left.+\frac{16\tau_0^4A_1c}{\tau_0^2\underline{\tau}^2}\right\}\label{rho2}\\
     \vartheta_1=&\max\left\{\frac{4\overline{\tau}^2b}{\tau_0^2},2b\overline{\tau}^2\left(\frac{30\overline{\tau}^2M_{v2_0}}{\beta_3}\right.\right.\nonumber\\
&\left.\left.+\frac{30\overline{\tau}^2(D_4M_{v1}+M_{v2})+M_{v1}}{\beta_1\beta_3}\right)\right\}\label{vartheta_1}\\
     \vartheta_2=&\max\left\{14b\overline{\tau}^2A_1(1+30D_4\overline{\tau}^2),\frac{120\overline{\tau}^4A_1b}{\tau_0^2},\frac{120\overline{\tau}^4D_2b}{\tau_0^2},\right.\nonumber\\
&\left.120b\overline{\tau}^4(D_7+\lambda^2D_1+D_2A_1+A_1A_3)\right\}\label{vartheta_2}.
 \end{align}
 The following holds 
     \begin{align}\label{inequ-V2bis}
	\mathcal{L}_iV_3(t)<&-(\theta_1-\rho_1|\tau_0-\tau_i|^2-\vartheta_1\sum_{j=1}^{r}|q_{ij}(t)|\nonumber\\
 &\cdot(|\tau_0-\tau_i|^2+|\tau_0-\tau_j|^2)) V_3(t)\nonumber\\
 &-(\theta_2-\rho_2|\tau_0-\tau_i|^2-\vartheta_2\sum_{j=1}^{r}|q_{ij}(t)|\nonumber\\
 &\cdot(|\tau_0-\tau_i|^2+|\tau_0-\tau_j|^2))V_3(0).
	\end{align}
	Combining \eqref{operator-L-simple}, \eqref{inequ-V2bis} and the property of transition rate $q_{ij}(t)$, we arrive at
\begin{align}\label{inequ-LV3}
&\mathcal{L}V_3(t)=\sum_{j=1}^{r}P_{ij}(0,t)\mathcal{L}_jV_3(t)\nonumber\\
\le&-\left(\theta_1-(\rho_1+r\vartheta_1q^\star)\mathbb{E}_{[0,\tau(0)]}\left(|\tau_0-\tau(t)|^2\right)\right.\nonumber\\
&\left.-\vartheta_1\sum_{l=1}^{r}|\tau_0-\tau_l|^2\sum_{j=1}^{r}P_{ij}(0,t)|q_{jl}(t)|\right) V_3(t)\nonumber\\
&-\left(\theta_2-(\rho_2+r\vartheta_2q^\star)\mathbb{E}_{[0,\tau(0)]}\left(|\tau_0-\tau(t)|^2\right)\right.\nonumber\\
&\left.-\vartheta_2\sum_{l=1}^{r}|\tau_0-\tau_l|^2\sum_{j=1}^{r}P_{ij}(0,t)|q_{jl}(t)|\right)V_3(0)\nonumber\\
\le&-\left(\theta_1-(\rho_1+(r+1)\vartheta_1q^\star)\mathbb{E}_{[0,\tau(0)]}\left(|\tau_0-\tau(t)|^2\right)\right.\nonumber\\
&\left.-\vartheta_1\sum_{l=1}^{r}|\tau_0-\tau_l|^2\sum_{j=1,j\ne l}^{r}P_{ij}(0,t)q_{jl}(t)\right) V_3(t)\nonumber\\
&-\left(\theta_2-(\rho_2+(r+1)\vartheta_2q^\star)\mathbb{E}_{[0,\tau(0)]}\left(|\tau_0-\tau(t)|^2\right)\right.\nonumber\\
&\left.-\vartheta_2\sum_{l=1}^{r}|\tau_0-\tau_l|^2\sum_{j=1,j\ne l}^{r}P_{ij}(0,t)q_{jl}(t)\right)V_3(0)\nonumber\\
\le&-\left(\theta_1-(\rho_1+(r+2)\vartheta_1q^\star)\mathbb{E}_{[0,\tau(0)]}\left(|\tau_0-\tau(t)|^2\right)\right.\nonumber\\
&\left.-\vartheta_1\sum_{l=1}^{r}|\tau_0-\tau_l|^2\frac{\partial P_{il}(0,t)}{\partial t}\right) V_3(t)\nonumber\\
&-\left(\theta_2-(\rho_2+(r+2)\vartheta_2q^\star)\mathbb{E}_{[0,\tau(0)]}\left(|\tau_0-\tau(t)|^2\right)\right.\nonumber\\
&\left.-\vartheta_2\sum_{l=1}^{r}|\tau_0-\tau_l|^2\frac{\partial P_{il}(0,t)}{\partial t}\right)V_3(0).
\end{align}
Denote 
\begin{align}
  G(t)=&\sum_{l=1}^{r}|\tau_0-\tau_l|^2\frac{\partial P_{il}(0,t)}{\partial t} \\
  \phi(t)=&\theta-(\rho+(r+2)\vartheta q^\star)\mathbb{E}_{[0,\tau(0)]}\left(|\tau_0-\tau(t)|^2\right)\nonumber\\
  &-\vartheta G(t)
\end{align}
where $\theta=\min\{\theta_k\}$, $\rho=\max\{ \rho_k\}$, $\vartheta=\max \{\vartheta_k\}$, $k=1,2$. Then,
we can rewrite \eqref{inequ-LV3} as
\begin{align}\label{L-2}
\mathcal{L}V_3(t)\le&-\phi(t)( V_3(t)+ V_3(0)).
\end{align}
Then, let
\begin{align}
F(t)=&e^{\int_0^t\phi(y)dy}(V_3(t)+V_3(0)).
\end{align}
Based on \eqref{L-2}, the following estimate holds
\begin{align}
\mathbb{E}_{[0,(V_3(0),\tau(0))]}(\mathcal{L}V_3(t)+\phi(t)V_3(t)+V_3(0))\le0
\end{align}
from which we can infer that the expectation of the infinitesimal operator $\mathcal{L}$ of $F(t)$ satisfies
\begin{align}
\mathbb{E}_{[0,(F(0),\tau(0))]}(\mathcal{L}F(t))\le&0.
\end{align}
Employing  the Dynkin's formula \cite{dynkin1965markov}, it can also get
\begin{align}\label{Dykin-equat}
\mathbb{E}_{[0,\tau(0)]}(F(t))-F(0)=&\mathbb{E}_{[0,(F(0),\tau(0))]}(\mathcal{L}F(t))\le0.
\end{align}
Since
\begin{align}
\int_0^tG(y)dy=&\int_0^t\sum_{l=1}^{r}|\tau_0-\tau_l|^2\frac{\partial P_{il}(0,y)}{\partial y}dy\nonumber\\
\le&\mathbb{E}_{[0,\tau(0)]}\left(|\tau_0-\tau(t)|^2\right)
\end{align}
we can get
\begin{align}
&\mathbb{E}_{[0,(F(0),\tau(0))]}(F(t))\nonumber\\
=&\mathbb{E}_{[0,(F(0),\tau(0))]}\left(e^{\int_0^t\phi(y)dy}(V_3(t)+V_3(0))\right)\nonumber\\
=&\mathbb{E}_{[0,(F(0),\tau(0))]}\left(\exp\left(\int_0^t\left(\theta-(\rho+(r+2)\vartheta q^\star)\right.\right.\right.\nonumber\\
&\left.\cdot\mathbb{E}_{[0,\tau(0)]}\{(\tau_0-\tau(y))^2\}-\vartheta G(y)\right)dy)\nonumber\\
&\left.\cdot( V_3(t)+V_3(0))\right)\nonumber\\
\ge&\mathbb{E}_{[0,(F(0),\tau(0))]}\left(e^{\int_0^t\left(\theta-(\rho+(r+2)\vartheta q^\star)\mathbb{E}_{[0,\tau(0)]}\{(\tau_0-\tau(y))^2\}\right)dy}\right.\nonumber\\
&\left.\cdot e^{-\vartheta\mathbb{E}_{[0,\tau(0)]}\left(|\tau_0-\tau(t)|^2\right)}(V_3(t)+V_3(0))\right)\nonumber\\
\ge&\mathbb{E}_{[0,(F(0),\tau(0))]}\left(e^{\int_0^t\left(\theta-(\rho+(r+2)\vartheta q^\star)\mathbb{E}_{[0,\tau(0)]}\{(\tau_0-\tau(y))^2\}\right)dy}\right.\nonumber\\
&\left.\cdot e^{-\vartheta \mathbb{E}_{[0,\tau(0)]}\left(|\tau_0-\tau(t)|^2\right)}V_3(t)\right).
\end{align}
Let $\mathbb{E}_{[0,\tau(0)]}\left(|\tau_0-\tau(t)|^2\right)\le \varepsilon$,   
\begin{align}
    \varepsilon (\lambda, q^\star)=\frac{\theta(\lambda)}{2\left[\rho(\lambda)+q^\star(r+2)\vartheta (\lambda)\right]}.\label{value-epsilon}
\end{align}
Indeed, $\lambda$ is a constituent variable of  $\varepsilon$ as the constants $A_i$ and $D_i$ depend on the gain kernels, which is essentially parameterized by the reactivity coefficient $\lambda$ according to Lemma \ref{inequ-k0} and Lemma \ref{inequ-k1-gamma1}.
Using \eqref{Dykin-equat} and \eqref{value-epsilon}, the following estimate holds
\begin{align}
&\mathbb{E}_{[0,(F(0),\tau(0))]}\left(e^{\frac{\theta}{2}t-\vartheta\varepsilon(\lambda, q^\star)}V_3(t)\right)\le F(0)= 2V_3(0).
\end{align}
 Finally, we get 
\begin{align}
\mathbb{E}_{[0,(V_3(0),\tau(0))]}(V_3(t))\le 2e^{\vartheta\varepsilon(\lambda, q^\star) } e^{-\frac{\theta}{2}t}V_3(0).
\end{align}

	Based on the equivalence between the Lyapunov functions $V_2(t)$ and $V_3(t)$,  a similar result can be established for $V_2(t)$, which concludes the proof of Proposition \ref{stbility-tsrget-m}.
\end{Proof}

According to Proposition \ref{prop2}, the norms of the target system \eqref{target-w-main}--\eqref{h-x-t} are equivalent to that of the cascade system \eqref{cascade-u-main}--\eqref{cascade-initial-tilde-v0}. Hence, combining  Proposition \ref{prop2}  and \ref{stbility-tsrget-m}, one can prove exponential stability of the system \eqref{cascade-u-main}--\eqref{cascade-initial-tilde-v0} as claimed in Theorem \ref{stability-cascade}, so is the original system \eqref{u-main}--\eqref{initial-u}.

\begin{Remark}\rm
    From \eqref{value-epsilon}, one can observe that $\varepsilon$ is a decreasing function of    $q^\star$, the maximum value of the transition rate  $q_{ij}(t)$: high values of  $q^\star$ lead to a closed-loop system that is less robust to stochastic error disturbances. The parameters $\theta_k$, $\rho_k$, and $\vartheta_k$ depend on the kernel functions $p(x,y),\ \gamma(x,y)$, and $k(x,y)$, as well as the upper and lower bounds of the stochastic delay  $\tau(t)$. As a result, the parameter $\varepsilon$ is directly impacted by the parameter $\lambda$, which is the exclusive parameter defining the gain kernels. A higher value of $\lambda$ leads to increased instability of the open-loop plant, necessitating a higher controller gain to restore stability to the disturbed closed-loop system.
\end{Remark}
\section{Simulation} \label{4}
We provide numerical examples to illustrate the effectiveness and examine the limitations of the delay-compensation controller in \cite{Krstic2009control} for system \eqref{cascade-u-main}--\eqref{cascade-initial-tilde-v0} when the delay is subject to stochastic error disturbances. 
We set $\lambda=11$, define the states of  
the stochastic delay $(\tau_1,\tau_2,\tau_3,\tau_4,\tau_5)=(0.3,0.4,0.5,0.6,0.7)$, the initial transition probabilities $(0.2,0.2,0.2,0.2,0.2)$, and 
the infinitesimal generator matrix  $Q=(q_{ij})_{r\times r}$
	\begin{align} \label{Q}      
	Q=
	\left(                 
	\begin{array}{ccccc}   
	-5 & 1 & 1 & 1 & 2\\  
	2 & -3 & 0.4 & 0.5 & 0.1\\  
	0.1 & 0.1 & -1 & 0.4 & 0.4\\
	2 & 1 & 0.5 & -4.5 & 1\\
	0.1 & 0.1 & 1 & 0.4 & -1.6\\
	\end{array}
	\right)              
	\end{align}
which reveals the instantaneous rates of transition between the stochastic delay values $\tau_i, i=1,2,\cdots, 5$. Note that the diagonal elements of $Q$ are chosen such that the row sums are zero to satisfy \eqref{10}.
\begin{figure}[t]
	\centering
	\subfloat[]{\label{delay1}\epsfig{figure=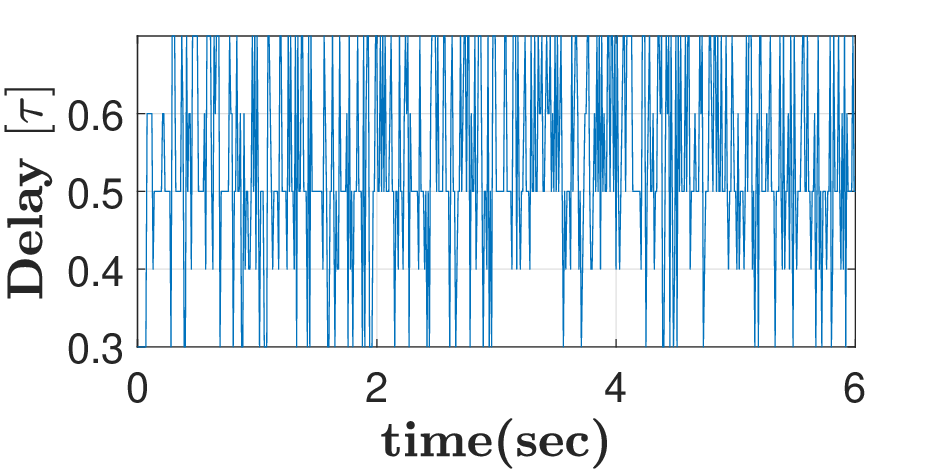,width=0.24\textwidth}}
	\subfloat[]{\label{delay2}\epsfig{figure=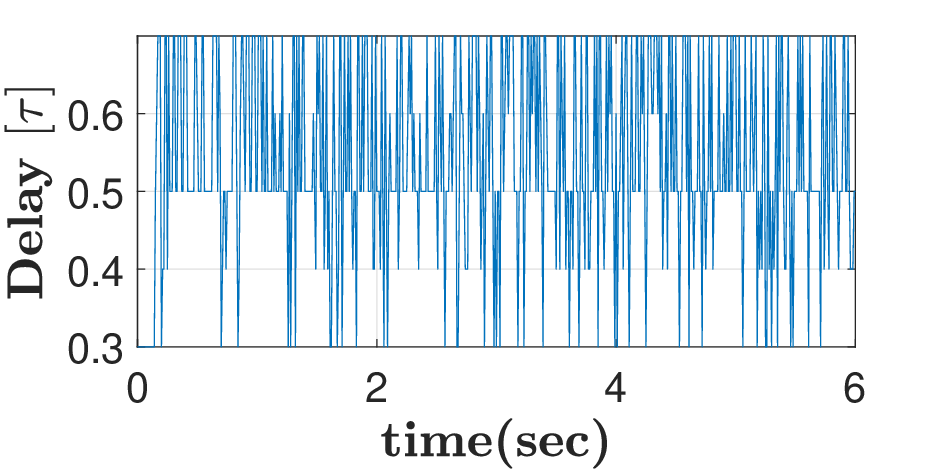,width=0.24\textwidth}}
	\caption{\small The realizations of stochastic delays. (a) Used in the first simulation. (b) Used in the second simulation.}\label{fig:delay}
\end{figure}
\begin{figure}[t]
	\centering
	\subfloat[]{\label{fig:stable-state}\epsfig{figure=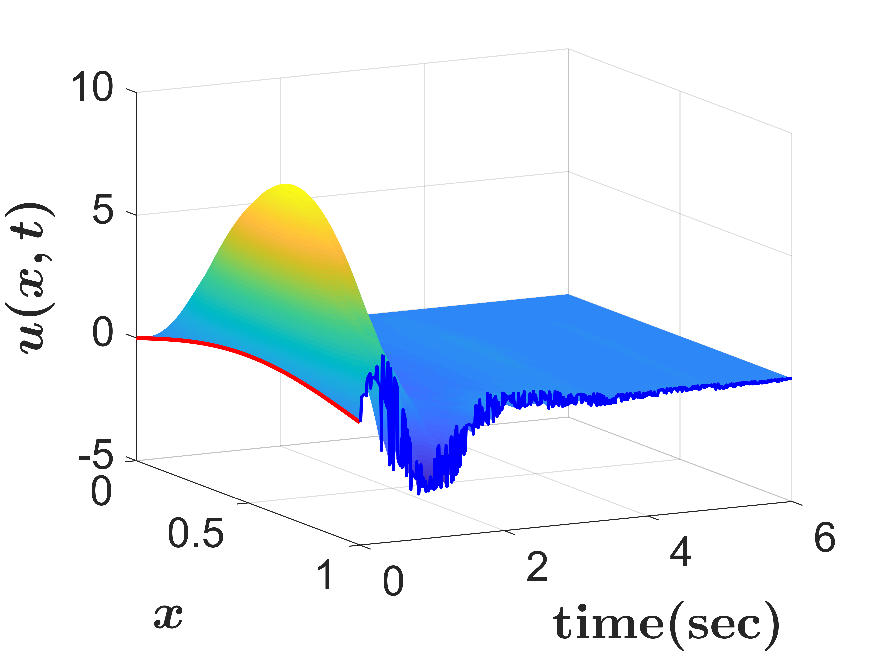,width=0.24\textwidth}}
	\subfloat[]{\label{fig:stable-norm}\epsfig{figure=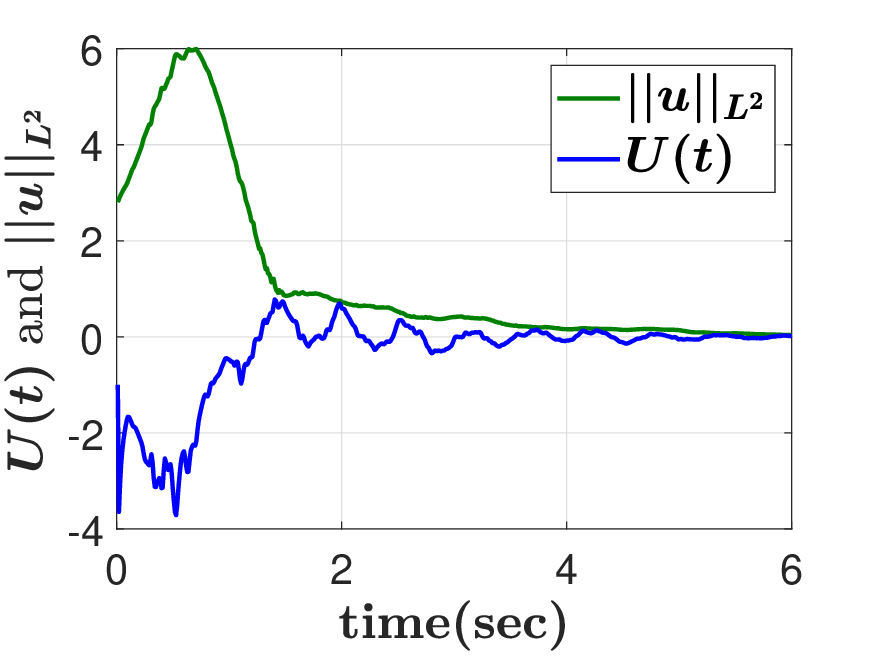,width=0.24\textwidth}}\\
 \subfloat[]{\label{fig:unstable-state}\epsfig{figure=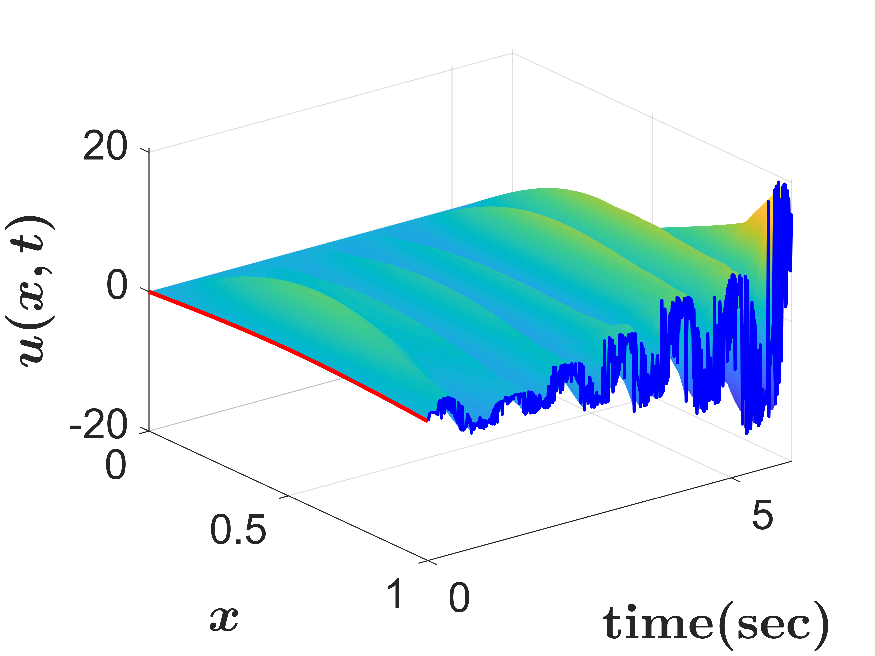,width=0.24\textwidth}}
	\subfloat[]{\label{fig:unstable-norm}\epsfig{figure=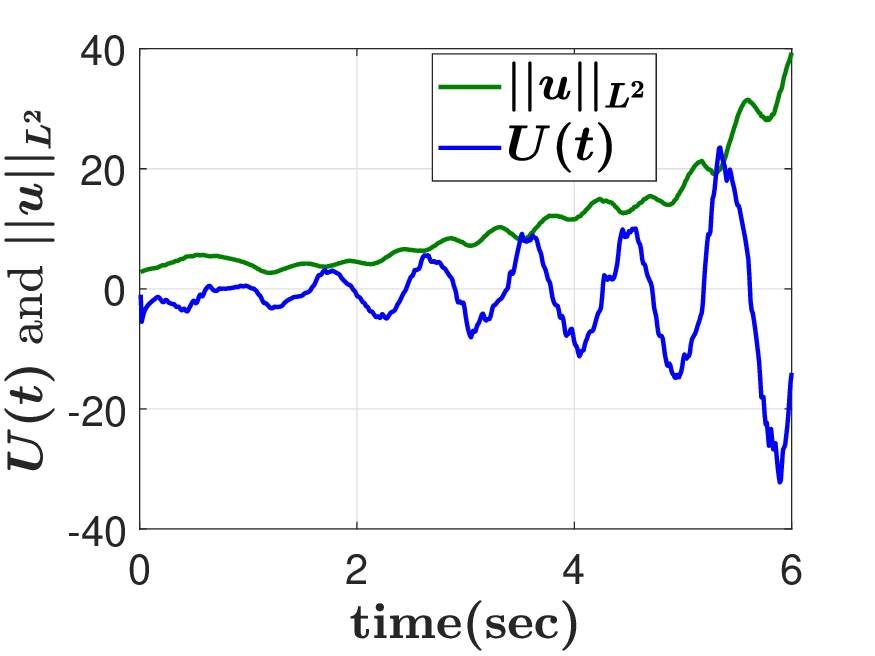,width=0.24\textwidth}}
	\caption{\small The simulation with $\lambda=11$ and (a), (b) the first simulation with $\tau_0=0.5$; (c), (d) the second simulation with $\tau_0=0.7$: (a), (c) the dynamics of the state  $u(x,t)$ under the controller $U(t)$;  (b), (d) the  $L^2$-norm of the state $u(x,t)$ under the controller $U(t)$ and the time-evolution of the input signal $U(t)$.}\label{fig:first}
\end{figure}

\begin{figure}[t]
	\centering
	\subfloat[]{\label{fig:stable-state-100}\epsfig{figure=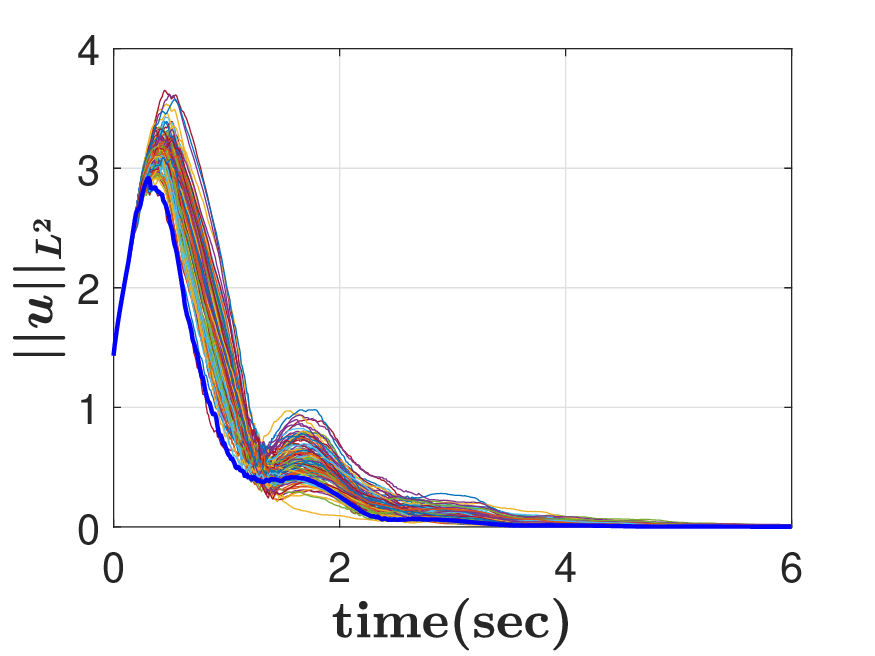,width=0.24\textwidth}}
	\subfloat[]{\label{fig:stable-input-100}\epsfig{figure=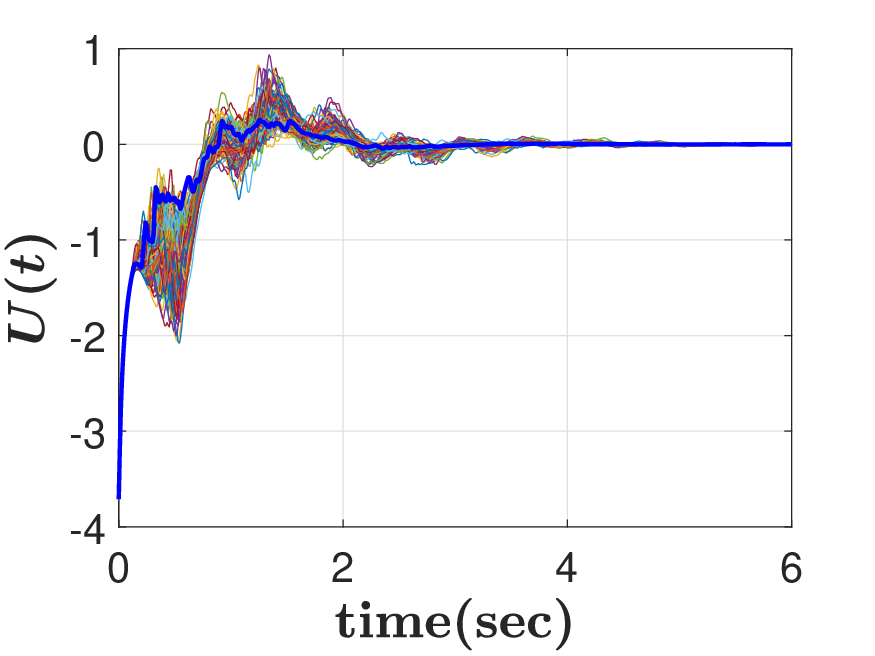,width=0.24\textwidth}}\\
 \subfloat[]{\label{fig:unstable-state-100}\epsfig{figure=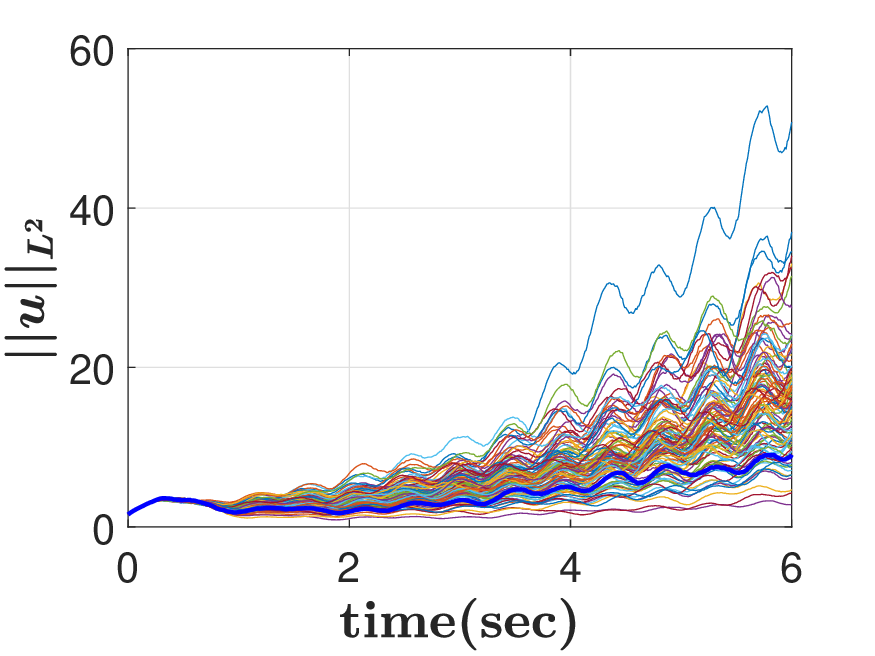,width=0.24\textwidth}}
	\subfloat[]{\label{fig:unstable-input-100}\epsfig{figure=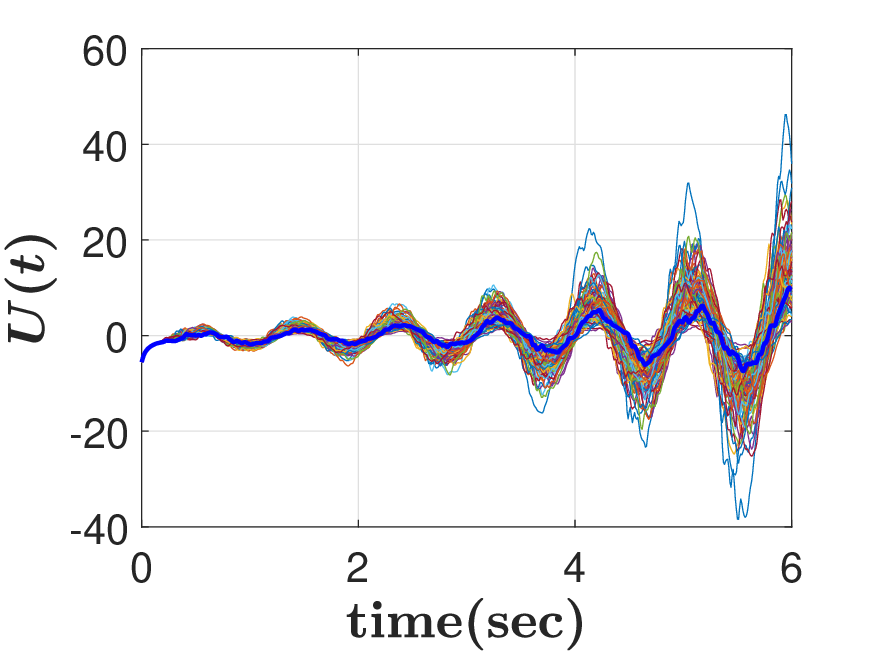,width=0.24\textwidth}}
	\caption{\small The simulation with $\lambda=11$ and (a), (b) with $\tau_0=0.5$; (b), (d) with $\tau_0=0.7$ under Monte Carlo method (100 trials): (a), (c) the $L^{2}$-norm of state $u(x,t)$;  (b), (d) the time-evolution of the input signal $U(t)$.}\label{fig:MC-stable}
\end{figure}

\begin{figure}[t]
	\centering
	\subfloat[]{\label{fig:la-stable-state}\epsfig{figure=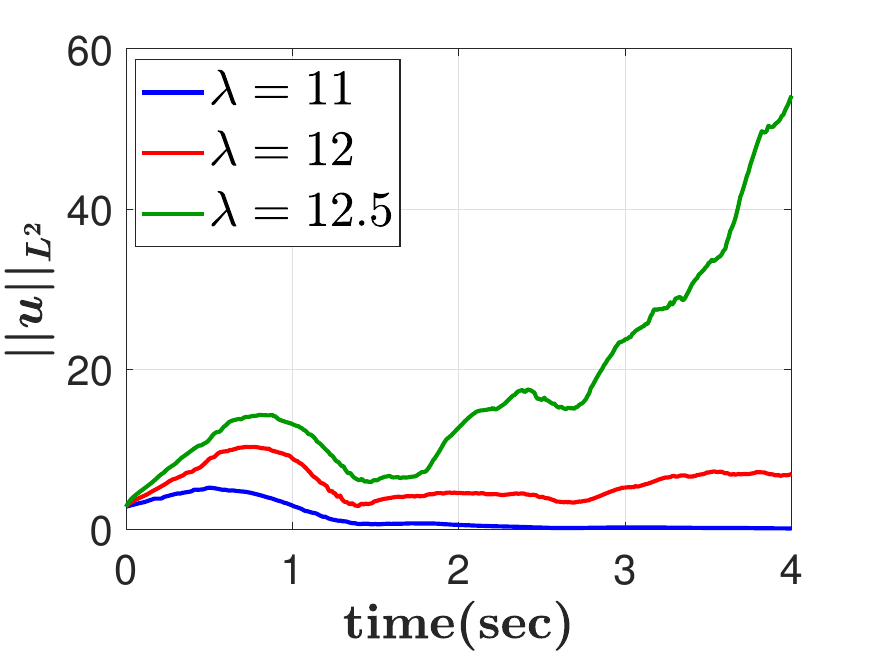,width=0.24\textwidth}}
 \subfloat[]{\label{fig:la-stable-norm}\epsfig{figure=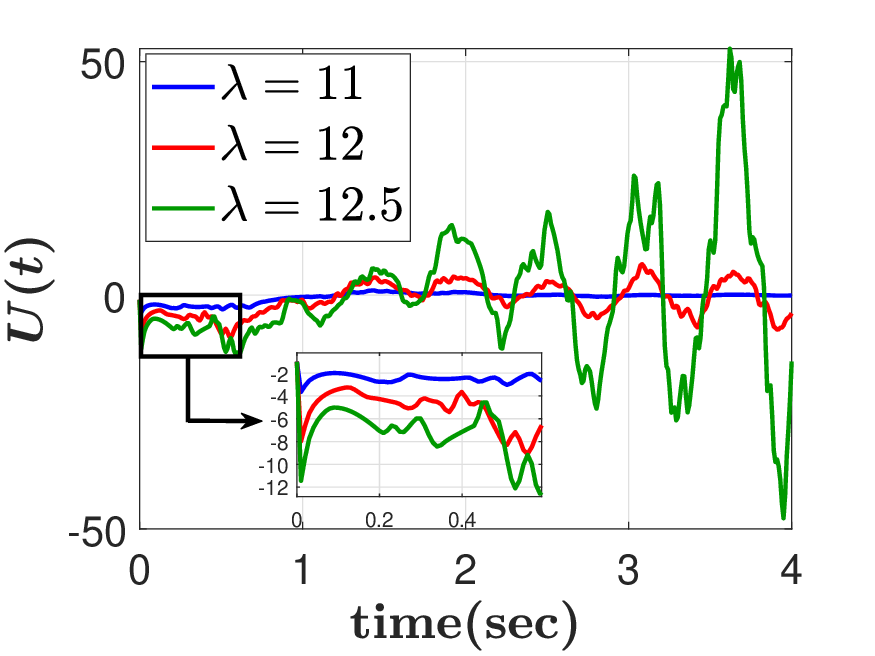,width=0.24\textwidth}}
	\caption{\small The simulation under the transition rate matrix $Q$ defined in \eqref{Q}, $\tau(t)\in[0.3,0.4,0.5,0.6,0.7]$ and $\tau_0=0.5$: (a) the $L^2$-norm of the system state $u(x,t)$;  
	 (b) the time evolution of the input signal $U(t)$.}\label{fig:lambda}
\end{figure}

\begin{figure}[t]
	\centering
	\subfloat[]{\label{fig:q-stable-state}\epsfig{figure=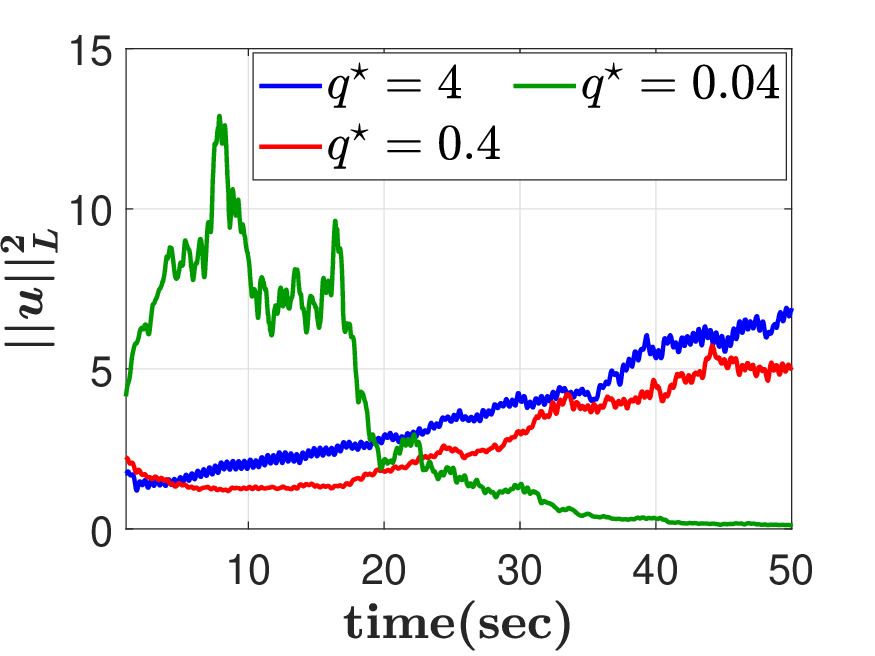,width=0.43\textwidth}}\\ \vspace{-0.1in}
 \subfloat[]{\label{fig:q-stable-norm}\epsfig{figure=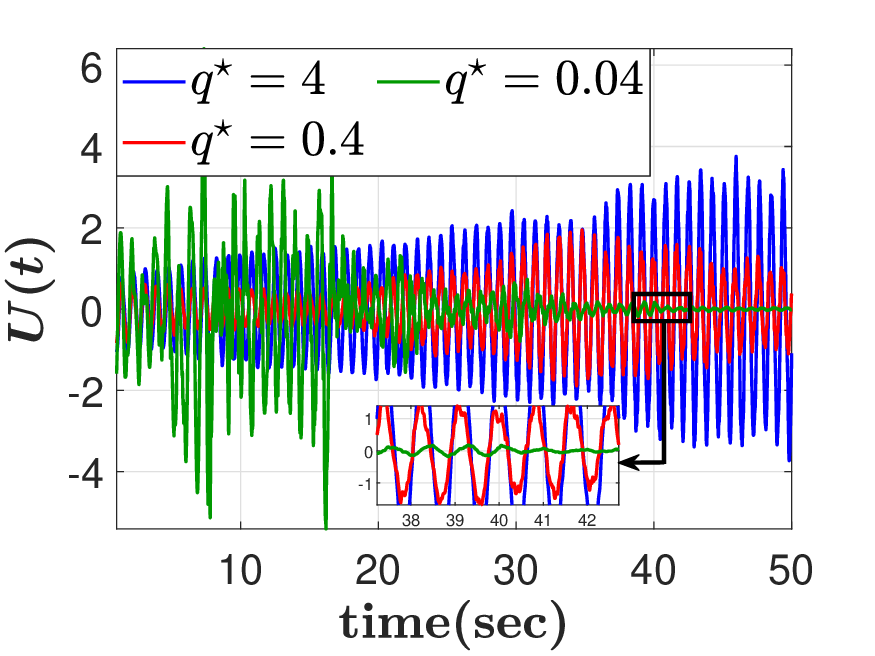,width=0.43\textwidth}}
	\caption{\small The simulation under $\lambda=11$, $\tau(t)\in[0.3,0.4,0.5,0.6,0.7]$ and $\tau_0=0.665$: (a) the $L^2$-norm of the system state $u(x,t)$;  
	 (b) the time evolution of the input signal $U(t)$.}\label{fig:qstar}
\end{figure}
The initial conditions are chosen as $u_0=\sin(\pi x)$, $\hat{v}_0=\cos(\pi x)$, $\bm{\tilde{v}}_0=(0.5\sin(\pi x),\sin(\pi x),1.5\sin(\pi x),2\sin(\pi x)$, $2.5\sin(\pi x))$. 
We discretize the system \eqref{cascade-u-main}--\eqref{cascade-initial-tilde-v0} by the finite difference method using the Crank-Nicolson, the spatial step is $\Delta x=0.01$ and the time step is $\Delta t=0.01s$ and use the Monte Carlo method \cite{ali2017multilevel} to simulate the closed-loop system with stochastic delay error disturbance. 

Two stochastic delay distributions are considered in the subsequent  simulation results:
\begin{itemize}
\item For the first scenario illustrated by  Fig. \ref{fig:delay} \subref{delay1}, with $\tau_0=0.5$ and an initial expectation set to $\mathbb{E}\left(|\tau_0-\tau(t)|^2\right)=0.02$ at $t=0$, the system exhibits mean-square $L^2$-exponential stability under the nominal control law \eqref{delay-controller}, as shown in Fig. \ref{fig:first} \subref{fig:stable-state} and \subref{fig:stable-norm}, respectively.
\item Conversely, the second scenario presented in Fig. \ref{fig:delay} \subref{delay2}, with $\tau_0=0.7$ and an initial expectation set to $\mathbb{E}\left(|\tau_0-\tau(t)|^2\right)=0.06$ at $t=0$, the closed-loop system exhibits instability  as depicted in Fig. \ref{fig:first} \subref{fig:unstable-state} and \subref{fig:unstable-norm}. Both the plant's state and the nominal control law \eqref{delay-controller}  diverge.  
\end{itemize}
These results confirm the theoretical findings as the expectation of the difference between stochastic delay and nominal delay must be reasonably small to guarantee the mean-square exponential stability of the system.  Further simulation results are generated to support our statement  using  a Monte
Carlo simulation with 100 trials in  Fig. \ref{fig:MC-stable}. As expected, the resulting simulations consolidate stability when the stochastic delay evolves in a sufficiently small interval  $\varepsilon$ (see. Fig. \ref{fig:stable-state-100} and \ref{fig:stable-input-100}) and instability otherwise (see. Fig. \ref{fig:unstable-state-100} and \ref{fig:unstable-input-100}). 

The influence of parameters $\lambda$ and $q^\star$ on the stability of the system \eqref{cascade-u-main}--\eqref{cascade-initial-tilde-v0} under the control law  \eqref{delay-controller} is depicted in Fig. \ref{fig:lambda} and Fig. \ref{fig:qstar}. Under consistent stochastic delay's variations, higher values of parameter $\lambda$ correspond to reduced control effectiveness. Similarly, when the parameter $q^\star$ increases, it diminishes the controller's capacity to effectively regulate the perturbed system: a smaller $\varepsilon$ decreases the system's tolerance to stochastic deviations from the nominal delay $\tau_0.$

\section{Conclusion}\label{5}
This paper is a study of  the robustness of the PDE backstepping delay-compensated controller for reaction-diffusion PDEs \cite{Krstic2009control} when stochastic disturbances act on the delay parameter. As compared to the previous works \cite{kong2022prediction1} and \cite{kong2022prediction2}, which compensate for stochastic delay for the systems modeled by  ODEs, our work extends the notion of robustness of predictor-feedback control laws  to infinite-dimensional class of systems when the delay is subject to Markov jumps perturbations. By establishing the well-posedness of the closed-loop system and conducting a stability analysis, our proof  reaches high levels of complexity. The maximum  deviation induced by the stochastic uncertainties    relative to the nominal delay is thoroughly estimated by assuming a known upper limit on the transition rate. 
Future developments will examine stochastic delay compensation when only boundary measurements are available, since measuring a spatially distributed state in real-time might be expensive.  %


\bibliographystyle{elsarticle-num}
\bibliography{reference}

\begin{appendices}
	\section{The proof of Lemma \ref{inequ-k0}}\label{proof-k0}
	Here, we give the proof of Lemma \ref{inequ-k0}.
	\begin{Proof}\rm
		Using the Cauchy-Schwarz inequality, the Fourier series of $p_y(1,y)$ and Lemma \ref{inequ-p}, one gets
		\begin{align*}
		\Vert k(x,0)\Vert^2=&\int_0^1\left(\sum_{n=1}^{\infty}(-1)^nn\pi e^{\tau_0(\lambda-n^2\pi^2)x}p_{1n}\right)^2dx\\
		\le&\sum_{n=1}^{\infty}(n\pi)^2p^2_{1n} \sum_{n=1}^{\infty}\int_0^1e^{2\tau_0(\lambda-n^2\pi^2)x}dx\\
		\le&\sum_{n=1}^{\infty}\frac{(e^{2\tau_0(\lambda-n^2\pi^2)}-1)}{\tau_0(\lambda-n^2\pi^2)}\Vert p_y(1,y)\Vert^2_{L^2}\le A_1.
		\end{align*}
		Similarly, one can deduce that
		\begin{align*}
		\Vert k_x(x,0)\Vert^2\le&\sum_{n=1}^{\infty}\frac{(e^{2\tau_0(\lambda-n^2\pi^2)}-1)}{\tau_0(\lambda-n^2\pi^2)}\Vert p_{yyy}(1,y)\Vert^2_{L^2}\le A_2.
		\end{align*}
		For the third inequality, utilizing the discrete form of the Cauchy-Schwarz inequality, the Fourier series of $p_{yy}(1,y)$ and Lemma \ref{inequ-p}, we get
		\begin{align*}
		k^2(0,0)=&\left(\sum_{n=1}^{\infty}(-1)^nn\pi p_{1n}\right)^2
		\le\sum_{n=1}^{\infty}\frac{1}{(n\pi)^2}\sum_{n=1}^{\infty}(n\pi)^4p_{1n}^2\\
		=&2\sum_{n=1}^{\infty}\frac{1}{(n\pi)^2}\Vert p_{yy}(1,y)\Vert^2_{L^2}\le A_3.
		\end{align*}
		The other inequalities i.e., $\Vert\eta(1,x)\Vert^2\le A_4$, $\Vert l_x(1,x)\Vert^2\le A_5$, $\Vert\eta_{xx}(1,x)\Vert^2\le A_6$, $|l(x,x)|^2\le A_7$, $\Vert l_{xx}(1,x)\Vert^2\le A_8$ can be proven by a similar approach, which concludes the proof.
	\end{Proof}
	\section{Some supports for Theorem \ref{well-t}}\label{A-Inner}
	\begin{Lemma}\label{inequs-inner}\rm
		Consider the function $d(x)\in L^2[0,1]$ defined in \eqref{d}, the following estimates hold
		\begin{align}
		&d^2(1)\le\left(5\tau^2_0A_1+\frac{5A_4}{3}\right)h_i^2(0)+5A_4\left(\int_{0}^{1}f^2dy\right.\nonumber\\
		&\left.+\frac{1}{3}g^2(0)\right)+5\tau^2_0A_5\int_{0}^{1}g^2dy\label{inequ-d1}\\
		&(d'(1))^2\le\frac{8\tau^4_0A_1}{\underline{\tau}^2}(h'_i(0))^2+8\tau^2_0A_8\int_{0}^{1}g^2dy\nonumber\\
		&+8\tau^4_0A_1\frac{(\max\delta_{\tau})^2}{\tau_0^2\underline{\tau}^2}d^2(0)+8\left(\tau^2_0A_2+\frac{1}{3}A_6+\tau^4_0A_1A_7\right)\nonumber\\
		&\cdot h^2_i(0)+8A_6\left(\int_{0}^{1}f^2dy+\frac{1}{3}g^2(0)\right).\label{inequ-dx1}
		\end{align}
	\end{Lemma}
	\begin{Proof} \rm
		Applying  the Cauchy-Schwarz inequality  to \eqref{d1} and using Lemma \ref{inequ-k0}, estimate  \eqref{inequ-d1} is deduced as follows 
		\begin{align*}
		d^2(1)=&\left(-\tau_0k(1,0)e^T_i\bm{h}(0)+\int_{0}^{1}\eta_x(1,y)(f(y)+yg(0,t)\right.\\
		&\left.+ye_i^T\bm{h}(0))dy+\tau_0\int_{0}^{1}l_x(1,y)g(y)dy\right)^2\\
		\le&\left(5\tau^2_0A_1+\frac{5A_4}{3}\right)h_i^2(0,t)+5A_4\left(\int_{0}^{1}f^2dy\right.\\
		&\left.+\frac{1}{3}g^2(0)\right)+5\tau^2_0A_5\int_{0}^{1}g^2dy.
		\end{align*}	
		Then, recalling \eqref{stable-z-main} and \eqref{stable-z-boundary}, one gets the following relations
		\begin{align*}
		g'(1)=&-\tau_0k(1,0)e_i^T\bm{h}(0)\\
		g''(1)=&-\tau^2_0k(1,0)e_i^T\Lambda^{-1}_{\tau}\left(\bm{h}'(0)-\Sigma_{\tau}d(0)\right)\nonumber\\
		&-\tau_0k_x(1,0)e_i^T\bm{h}(0).
		\end{align*}
		Combining the \eqref{hat-h} and the definition of $d(x)$, the following equality can be stated
  \begin{align}
      d(x)=&g'(x)+\int_0^1\eta_x(x,y)(f(y)+yg(0)+ye_i^T\bm{h}(0))dy\nonumber\\
      &+\tau_0l(x,x)g(x)+\tau_0\int_0^1l_x(x,y)g(y)dy
  \end{align}
  and with the help of Lemma \ref{inequ-k0}, one can infer that
		\begin{align*}
		&(d'(1))^2=\left(g''(1)+\tau_0l(1,1)g'(1)+\int_{0}^{1}\eta_{xx}(1,y)(f(y)\right.\\
		&+yg(0)+ye_i^T\bm{h}(0))dy\left.+\tau_0\int_{0}^{1}l_{xx}(1,y)g(y)dy\right)^2\\
		&\le\frac{8\tau^4_0A_1}{\underline{\tau}^2}(h'_i(0))^2+8\tau^2_0A_8\int_{0}^{1}g^2dy+8\left(\tau^2_0A_2+\frac{1}{3}A_6\right.\nonumber\\
		&\left.+\tau^4_0A_1A_7\right)h^2_i(0)+8\tau^4_0A_1\frac{(\max\delta_{\tau})^2}{\tau_0^2\underline{\tau}^2}d^2(0)\nonumber\\
		&+8A_6\left(\int_{0}^{1}f^2dy+\frac{1}{3}g^2(0)\right)
		\end{align*}
		to conclude the proof.
	\end{Proof}
	\section{Some supports for Proposition \ref{stbility-tsrget-m}}  \label{A-h}
	
	\begin{Lemma}\label{inequs-V}\rm
		Consider the function $\hat{h}(x,t)$ defined in \eqref{hat-h}, the following inequalities hold
		\begin{align}
		&\hat{h}^2(1,t)\le\left(5\tau^2_0A_1+\frac{5A_4}{3}\right)\tilde{v}_i^2(0,t)+5A_4\left(\int_{0}^{1}m^2(y,t)dy\right.\nonumber\\
		&\left.+\frac{1}{3}\hat{z}^2(0,t)\right)+5\tau^2_0A_5\int_{0}^{1}\hat{z}^2(y,t)dy\label{inequ-h1}\\
		&\hat{h}^2_x(1,t)\le\frac{8\tau^4_0A_1}{\underline{\tau}^2}(\partial_x\tilde{v}_i(0,t))^2+8\tau^4_0A_1\frac{(\tau_0-\tau_i)^2}{\tau_0^2\underline{\tau}^2}\nonumber\\
		&\cdot\hat{h}^2(0,t)+8\left(\tau^2_0A_2+\frac{1}{3}A_6+\tau^4_0A_1A_7\right)\tilde{v}^2_i(0,t)+8A_6\nonumber\\
		&\cdot\left(\int_{0}^{1}m^2(y,t)dy+\frac{1}{3}\hat{z}^2(0,t)\right)+8\tau^2_0A_8\int_{0}^{1}\hat{z}^2(y,t)dy.\label{inequ-hx1}
		\end{align}
	\end{Lemma}
	The proof of Lemma \ref{inequs-V} can be established by following the steps of the proof of Lemma \ref{inequs-inner}.

	\begin{Lemma}\rm\label{bound-V3-ij}
		Consider the Lyapunov function $V_3(t)$ given by \eqref{V-2}, the following inequality holds
		\begin{align}
		&\sum_{j=1}^{r}q_{ij}(t)(V_3(t,j)-V_3(t,i))\nonumber\\
		\le& 2b\overline{\tau}^2\sum_{j=1}^{r}|q_{ij}(t)|(|\tau_0-\tau_i|^2+|\tau_0-\tau_j|^2)\nonumber\\
  &\cdot\left(7A_1(1+30D_4\overline{\tau}^2)\hat{z}^2(0,t)+\frac{2}{\tau_0^2}\hat{h}^2(x,t)\right.\nonumber\\
  &+\left(\frac{30\overline{\tau}^2M_{v2_0}}{\beta_3}+\frac{30\overline{\tau}^2(D_4M_{v1}+M_{v2})+M_{v1}}{\beta_1\beta_3}\right)V_3(t)\nonumber\\
  &+30\overline{\tau}^2\left(\frac{2D_2}{\tau_0^2}(\partial_x\tilde{v}_i(0,t))^2+2(D_7+\lambda^2D_1+D_2A_1\right.\nonumber\\
&\left.\left.+A_1A_3)\tilde{v}_i^2(0,t)+\frac{2A_1}{\tau^2_0}\hat{z}^2_x(0,t)\right)\right).
		\end{align}
	\end{Lemma}
	\begin{Proof}\rm
		Based on \eqref{V-2}, the following is deduced
		\begin{align*}
		&\sum_{j=1}^{r}q_{ij}(t)(V_3(t,j)-V_3(t,i))\nonumber\\
		=&b\sum_{j=1}^{r}q_{ij}(t)\int_0^1(1+x)\left(\tau^2_j\tilde{v}_j^2(x,t)-\tau^2_i\tilde{v}_i^2(x,t)\right)dx\nonumber\\
  &+b\sum_{j=1}^{r}q_{ij}(t)\int_0^1(1+x)\left(\tau^2_j(\partial_x\tilde{v}_j(x,t))^2\right.\nonumber\\
  &\left.-\tau^2_i(\partial_x\tilde{v}_i(x,t))^2\right)dx.
		\end{align*}
		Combining the definition of $\bm{\tilde{v}}(x,t)$ in \eqref{error-control} leads to the following estimate
		\begin{align*}
		&\tau^2_j\tilde{v}_j^2(x,t)-\tau^2_i\tilde{v}_i^2(x,t)\nonumber\\
		=&\left(\tau_j\int_{t+\tau_0(x-1)}^{t+\tau_j(x-1)}\dot{U}(s)ds-\tau_i\int_{t+\tau_0(x-1)}^{t+\tau_i(x-1)}\dot{U}(s)ds\right)\nonumber\\
		&\cdot\left(\tau_j\int_{t+\tau_0(x-1)}^{t+\tau_j(x-1)}\dot{U}(s)ds+\tau_i\int_{t+\tau_0(x-1)}^{t+\tau_i(x-1)}\dot{U}(s)ds\right)\nonumber\\
  \le&2\overline{\tau}^2\left(|\tau_j-\tau_0|^2+|\tau_i-\tau_0|^2\right)\left(\max_{s\in[0,t]}\dot{U}(s)\right)^2.
		\end{align*}
  Then, using the definition of $\bm{\tilde{v}}(x,t)$ in \eqref{error-control} and \eqref{stable-ev-main}, we arrive at
  \begin{align*}
      &\tau^2_j\partial_x\tilde{v}_j^2(x,t)-\tau^2_i\partial_x\tilde{v}_i^2(x,t)\nonumber\\
      =&\left(\tau_j^2\int_{t+\tau_0(x-1)}^{t+\tau_j(x-1)}\ddot{U}(s)ds-\tau_i^2\int_{t+\tau_0(x-1)}^{t+\tau_i(x-1)}\ddot{U}(s)ds\right.\nonumber\\
      &\left.+\frac{(\tau_j-\tau_0)\tau_j}{\tau_0}\hat{h}(x,t)-\frac{(\tau_i-\tau_0)\tau_i}{\tau_0}\hat{h}(x,t)\right)\nonumber\\
		&\cdot\left(\tau_j^2\int_{t+\tau_0(x-1)}^{t+\tau_j(x-1)}\ddot{U}(s)ds+\tau_i^2\int_{t+\tau_0(x-1)}^{t+\tau_i(x-1)}\ddot{U}(s)ds\right.\nonumber\\
      &\left.+\frac{(\tau_j-\tau_0)\tau_j}{\tau_0}\hat{h}(x,t)+\frac{(\tau_i-\tau_0)\tau_i}{\tau_0}\hat{h}(x,t)\right)\nonumber\\
  \le&4\left(\overline{\tau}^4\left(|\tau_j-\tau_0|^2+|\tau_i-\tau_0|^2\right)\left(\max_{s\in[0,t]}\ddot{U}(s)\right)^2\right.\nonumber\\
  &\left.+\frac{(|\tau_j-\tau_0|^2+|\tau_i-\tau_0|^2)\overline{\tau}^2}{\tau_0^2}|\hat{h}(x,t)|^2\right)
  \end{align*}
		from which one obtains
		\begin{align*}
		&\sum_{j=1}^{r}q_{ij}(t)(V_3(t,j)-V_3(t,i))\nonumber\\
		\le& 2b\overline{\tau}^2\sum_{j=1}^{r}|q_{ij}(t)|(|\tau_j-\tau_0|^2+|\tau_i-\tau_0|^2)\left(\frac{2}{\tau_0^2}|\hat{h}(x,t)|^2\right.\nonumber\\
  &\left.+2\overline{\tau}^2\left(\max_{s\in[0,t]}\ddot{U}(s)\right)^2+\left(\max_{s\in[0,t]}\dot{U}(s)\right)^2\right).
		\end{align*}
		Next,   as $\delta(t)=e_i$,  with the help of \eqref{delay-controller}, \eqref{cascade-u-main}--\eqref{cascade-u-bnd1}, and \eqref{transformation-zv}, the Cauchy-Schwarz inequality, Lemma \ref{inequ-k0}, Lemma \ref{inequ-k1-gamma1}, Proposition \ref{prop2}, inequality \eqref{equ-V-2-3} and Poincar\'{e} inequality, i.e., $\hat{v}^2(0,t)\le\int_0^1 \hat{v}^2_x(x,t)dx$, $\tilde{v}_i(0,t)\le\int_0^1 (\partial_x\tilde{v}_i(x,t))^2dx$, the following holds
		\begin{align*}
		&(\dot{U}(t))^2\nonumber\\
		=&\left(\gamma_y(1,1)u(1,t)+\int_{0}^{1}(\gamma_{yy}(1,y)+\lambda\gamma(1,y))u(y,t)dy\right.\nonumber\\
		&+k(1,1)U(t)-k(1,0)\int_{0}^{1}\gamma(0,y)u(y,t)dy\nonumber\\
		&\left.-k(1,0)\hat{z}(0,t)-\int_{0}^{1}k_y(1,y)\hat{v}(y,t)dy\right)^2\nonumber\\
		\le&7\left(2\gamma^2_y(1,1)\left(\int_{0}^{1}(\partial_x\tilde{v}_i(x,t))^2dx+\int_{0}^{1}\hat{v}^2_x(x,t)dx\right)\right.\nonumber\\
		&+\left(\int_{0}^{1}\gamma^2_{yy}(1,y)dy+(\lambda^2+2k^2(1,1))\int_{0}^{1}\gamma^2(1,y)dy\right.\nonumber\\
		&\left.+ k^2(1,0)\int_0^1\gamma^2(0,y)dy\right)\int_{0}^{1}u^2(x,t)dx\nonumber\\
		&+\left(\int_{0}^{1}k^2_y(1,y)dy+2k^2(1,1)\tau^2_0\int_{0}^{1}k^2(1,y)dy\right)\nonumber\\
		&\left.\cdot\int_{0}^{1}\hat{v}^2(y,t)dy+k^2(1,0)\hat{z}^2(0,t)\right)\nonumber\\
		\le&7\left(2D_2\int_{0}^{1}(\partial_x\tilde{v}_i(x,t))^2dx+2D_2\int_{0}^{1}\hat{v}^2_x(x,t)dx\right.\nonumber\\
		&+\left(D_3+\lambda^2D_1+2D_1D_4+A_1M_1\right)\int_{0}^{1}u^2(x,t)dx\nonumber\\
		&\left.+(2\tau_0^2D_4^2+D_5)\int_{0}^{1}\hat{v}^2(y,t)dy+A_1\hat{z}^2(0,t)\right)\nonumber\\
  \le&\frac{M_{v1}}{\beta_1\beta_3} V_3(t)+7A_1\hat{z}^2(0,t)
		\end{align*}
        and
        \begin{align*}
		&\left(\ddot{U}(t)\right)^2\nonumber\\
		=&\left(-\gamma_y(1,0)\left(\hat{v}_t(0,t)+\partial_t\tilde{v}_i(0,t)\right)-(\gamma_{yyy}(1,1)\right.\nonumber\\
		&\left.+\lambda\gamma_y(1,1)\right)u(1,t)+\int_{0}^{1}(\gamma_{yyyy}(1,y)+2\lambda\gamma_{yy}(1,y)\nonumber\\
		&+\lambda^2\gamma(1,y))u(y,t)dy+k(1,1)\dot{U}(t)+k(1,0)\gamma_y(0,1)\nonumber\\
		&\cdot u(1,t)-k(1,0)\int_{0}^{1}(\gamma_{yy}(0,y)+\lambda\gamma(0,y))u(y,t)dy\nonumber\\
		&-k(1,0)\hat{z}_t(0,t)-\frac{1}{\tau_0}k_y(1,1)U(t)+\frac{1}{\tau_0}k_y(1,0)\hat{v}(0,t)\nonumber\\
  &\left.+\frac{1}{\tau_0}\int_0^1k_{yy}(1,y)\hat{v}(y,t)dy\right)^2\nonumber\\
		\le&15\left(\gamma_y^2(1,0)(\hat{v}^2_t(0,t)+(\partial_t\tilde{v}_i(0,t))^2)+(\gamma^2_{yyy}(1,1)\right.\nonumber\\
		&+\lambda^2\gamma^2_y(1,1)+k^2(1,0)\gamma^2_y(0,1)) u^2(1,t)\nonumber\\
		&+\left(\int_0^1\gamma^2_{yyyy}(1,y)dy+4\lambda^2\int_0^1\gamma^2_{yy}(1,y)dy+\lambda^4\right.\nonumber\\
  &\cdot\int_0^1\gamma^2(1,y)dy+k^2(1,0)\int_0^1\gamma^2_{yy}(0,y)dy+\lambda^2k^2(1,0)\nonumber\\
		&\left.\cdot\int_0^1\gamma^2(0,y)dy\right)\int_0^1u^2(y,t)dy+k^2(1,1)(\max\dot{U})^2\nonumber\\
		&+k^2(1,0)\hat{z}^2_t(0,t)+\frac{1}{\tau_0^2}k^2_y(1,1)U^2(t)+\frac{1}{\tau_0^2}k^2_y(1,0)\nonumber\\
		&\cdot\hat{v}^2(0,t)+\frac{1}{\tau_0^2}\int_0^1k^2_{yy}(1,y)dy\int_0^1\hat{v}^2(y,t)dy\nonumber\\
		\le&15\left(\left(\frac{D_2}{\tau^2_0}+\frac{2D_2(\tau_i-\tau_0)}{\tau^2_i\tau_0^2}\right)\hat{v}^2_x(0,t)+\frac{2D_2}{\tau_i^2}(\partial_x\tilde{v}_i(0,t))^2\right.\nonumber\\
		&+2\left(D_7+\lambda^2D_2+D_2A_1+\frac{D_5}{2\tau_0^2}\right)\hat{v}^2(0,t)+2\left(D_7\right.\nonumber\\
		&\left.+\lambda^2D_2+A_1D_2+A_1A_3\right)\tilde{v}_i^2(0,t)+\left(D_6+4\lambda^2D_3\right.\nonumber\\
		&\left.+\lambda^4D_1+A_1D_3+\lambda^2A_1D_1+\frac{2D_1D_5}{\tau_0^2}\right)\int_0^1u^2(y,t)dy\nonumber\\
		&+7A_1D_4\hat{z}^2(0,t)+D_4M_{v1}V_1(t)+\frac{2A_1}{\tau^2_0}\hat{z}^2_x(0,t)\nonumber\\
		&+\left(2D_4D_5+\frac{D_8}{\tau_0^2}\right)\int_0^1\hat{v}^2(y,t)dy\nonumber\\
		\le&15\left(\left(\frac{M_{v2_0}}{\beta_3}+\frac{D_4M_{v1}+M_{v2}}{\beta_1\beta_3}\right)V_3(t)\right.\nonumber\\
  &+\frac{2D_2}{\tau_0^2}(\partial_x\tilde{v}_i(0,t))^2+2\left(D_7+\lambda^2D_1+D_2A_1\right.\nonumber\\
		&\left.\left.+A_1A_3\right)\tilde{v}_i^2(0,t)+7A_1D_4\hat{z}^2(0,t)+\frac{2A_1}{\tau^2_0}\hat{z}^2_x(0,t)\right).
		\end{align*}
		The proof of Lemma \ref{bound-V3-ij} is complete.
	\end{Proof}

\end{appendices}

\end{document}